\documentclass{amsart}
\usepackage{graphicx}
\usepackage{amssymb}
\usepackage{epstopdf}
\usepackage{nicefrac}
\usepackage{enumerate}
\usepackage{hyperref}
\usepackage{float}
\usepackage{color}
\usepackage{pst-all}
\usepackage{tabularx}

\newtheorem{thm}{THEOREM}[section]

\newtheorem{cor}[thm]{COROLLARY}

\newtheorem{defn}[thm]{DEFINITION}
\newtheorem{ex}[thm]{EXAMPLE}

\newtheorem{lemma}[thm]{LEMMA}

\newtheorem{prop}[thm]{PROPOSITION}

\newtheorem{remark}[thm]{REMARK}

\newcommand{\cA}{{\mathcal A}}

\newcommand{\cC}{{\mathcal C}}
\newcommand{\cD}{{\mathcal D}}
\newcommand{\cG}{{\mathcal G}}
\newcommand{\cH}{{\mathcal H}}

\newcommand{\cN}{{\mathcal N}}
\newcommand{\cO}{{\mathcal O}}
\newcommand{\cP}{{\mathcal P}}
\newcommand{\cQ}{{\mathcal Q}}
\newcommand{\cS}{{\mathcal S}}
\newcommand{\cT}{{\mathcal T}}
\newcommand{\cW}{{\mathcal W}}

\newcommand{\ds}{\displaystyle}
\newcommand{\e}{{\varepsilon}} 
\newcommand{\G}{\Gamma}

\newcommand{\mZ}{{\mathbb Z}}
\newcommand{\whG}{{\widehat{G}}}
\newcommand{\oG}{{\overline{\Phi(G)}}}
\newcommand{\whg}{{\widehat{g}}}

\newcommand{\whx}{{\widehat{x}}}
\newcommand{\whcA}{{\widehat{\cA}}}
\newcommand{\whcN}{{\widehat{\cN}}}
\newcommand{\whTheta}{{\widehat{\Theta}}}
\newcommand{\whPhi}{{\widehat{\Phi}}}
\newcommand{\whE}{{\widehat{E}}}

\newcommand{\whphi}{\widehat{\phi}} %
\newcommand{\fG}{\mathfrak{G}} %

\input xy
\xyoption {all}

\begin{document}

\title{The discriminant invariant of Cantor group actions}

\begin{abstract}
In this work, we  investigate the dynamical and geometric properties of weak solenoids, as part of the development of a ``calculus of group chains'' associated to Cantor minimal actions. The       study of the properties  of group chains was initiated in the works of McCord \cite{McCord1965} and     Fokkink and Oversteegen   \cite{FO2002},   to   study the problem of determining which weak solenoids are homogeneous continua.  We develop an alternative condition for the homogeneity in terms of   the Ellis semigroup of the action, then investigate the relationship between   non-homogeneity  of a weak solenoid   and its discriminant   invariant, which we introduce in this work. A key part of our study is the construction of new examples that illustrate various subtle properties of group chains that correspond to geometric properties of non-homogeneous weak solenoids.
\end{abstract}

\author{Jessica Dyer}
\author{Steven Hurder}
 \author{Olga Lukina}
 \email{jesscdyer@gmail.com, hurder@uic.edu, lukina@math.uic.edu}
\address{Department of Mathematics, University of Illinois at Chicago, 322 SEO (m/c 249), 851 S. Morgan Street, Chicago, IL 60607-7045}
  
\thanks{Version date: May 2, 2016}



\maketitle


\section{Introduction}\label{sec-intro}

Let $X$ be a  Cantor set,   let $G$ be a finitely generated infinite group, and let $\Phi \colon G \times X \to X$ be an action by homeomorphisms. That is, there is a map $\Phi \colon G \to {\it Homeo}(X)$ which associates a homeomorphism $\Phi(g)$ of $X$ to each $g \in G$. We   adopt the short-cut notation $g \cdot x = \Phi(g)(x)$ when convenient. 
An action $(X,G,\Phi)$  is \emph{minimal} if for each $x \in X$, its orbit   $\cO(x) = \{g \cdot x \mid g \in G\}$ is dense in $X$. In this case, we say that    $(X,G,\Phi)$  is a \emph{Cantor minimal system}.

 Cantor minimal systems  $(X,G,\Phi)$ and $(Y,G, \Psi)$ are    \emph{topologically conjugate}, or just \emph{conjugate}, if there exists a homeomorphism $\tau\colon X \to Y$ such that $\tau(\Phi(g)(x)) = \Psi(g)(\tau(x))$. 
In this paper, we are concerned with the classification of  Cantor minimal systems up to topological conjugacy, and focus especially on invariants of the system   for the case when $G$ is non-abelian.

A topological space $\cS$ is \emph{homogeneous} if   for every $x,y \in \cS$ there is a homeomorphism $h \colon \cS  \to \cS$ such that $h(x) = y$.  One motivation for this work comes from the problem of classifying solenoids up to homeomorphism  as in \cite{ClarkHurder2013,CHL2013c}, and to understand their groups of self-homeomorphisms.

  Recall that     a \emph{weak solenoid} $\cS$ is an inverse limit of an increasing sequence of non-trivial finite-to-one   coverings $\pi_{\ell} \colon M_{\ell} \to M_0$ of a closed connected manifold $M_0$. The projection $f_0 \colon \cS \to M_0$ is a fiber bundle, whose typical fiber $F_0 = f_0^{-1}(x_0)$  is a Cantor set, for $x_0 \in M_0$. 
 The path lifting property for the finite covering maps $\pi_{\ell}$  induces the monodromy  action $\Phi \colon G_0 \times F_0 \to F_0$ of the 
   fundamental group $G_0 = \pi_1(M_0,x_0)$ on the inverse limit fiber $F_0$.
   Then $(F_0,G_0,\Phi)$ is  a Cantor minimal system.
   
A finite covering $\pi_{\ell} \colon M_{\ell} \to M_0$ is said to be \emph{regular}  if it is defined by a normal subgroup. That is, let $x_0 \in M_0$ be a chosen basepoint, and let $x_{\ell} \in M_{\ell}$ be a lift such that $\pi_{\ell}(x_{\ell}) = x_0$, then the covering is regular if the image of the induced map on fundamental groups, 
$(\pi_{\ell})_{\#} \colon \pi_1(M_{\ell}, x_{\ell}) \to \pi_1(M_0 , x_0)$, is a normal subgroup of $\pi_1(M_0 , x_0)$.
It is a standard fact (see \cite{Massey1967}) that a covering is regular if and only if the group   of deck transformations acts transitively on the   fibers of the covering.

If each   map $\pi_{\ell} \colon M_{\ell} \to M_0$ in the definition of a weak solenoid $\cS$ is a regular covering,   then we say that $f_0 \colon \cS \to M_0$ is a \emph{regular solenoid}.   The fiber $F_0$ of a regular solenoid is a Cantor group, and so there is a natural \emph{right} action of the Cantor group $F_0$ on $\cS$ which is transitive on fibers, and commutes with the left monodromy action of $G_0$ on the fibers.  McCord used this fact to show   in \cite{McCord1965}  that  a regular solenoid $\cS$ is  homogeneous.   Rogers and Tollefson \cite{RT1971b} subsequently  gave an example of a weak solenoid which is defined by a sequence of covering maps which are not regular coverings, but the inverse limit space $\cS$ is   still homogeneous. They posed the problem of determining, under which conditions is a weak solenoid $\mathcal{S}$   homogeneous?

On the other hand,   weak solenoids need not be homogeneous, and examples of such  were first given by Schori in \cite{Schori1966}, and later by   Rogers and Tollefson in \cite{RT1971b}.  Fokkink and Oversteegen   developed in \cite{FO2002} a criterion for a weak solenoid to be homogeneous, stated as Theorem~\ref{thm-weaklynormal} below, formulated in terms of  the properties of the nested group chain $\cG = \{ G_{i}  \mid i  \geq 0\}$ of subgroups of finite index in the group $G_0 = \pi_1(M_0, x_0)$, where $G_{i} \subset G_0$ is the image in $G_0$ of the fundamental group of the covering space  $M_{i}$.

The later work of Clark, Fokkink and Lukina \cite{CFL2014}  gave examples of weak solenoids for which the leaves (the path connected components) in the solenoid have different end structures,  and thus these solenoids cannot be homogeneous. The methods used in this paper were geometric in nature, and  showed that the geometry of the leaves in a weak solenoid are also obstacles to homogeneity.

The examples of Schori, Rogers and Tollefson, and Clark, Fokkink and Lukina, suggest a variety of   questions about the structure of weak solenoids.   In this work, we develop tools for their study.

We recall two properties of a Cantor minimal system which will be important for the following.    
First,  a   Cantor minimal system $(X,G,\Phi)$ is \emph{equicontinuous} with respect to a metric $d_X$ on $X$, if for all $\e >0$ there exists $\delta > 0$, such that for all $x , y \in X$ and $g \in G$ we have
  $$d_X(x,y) < \delta \qquad \Longrightarrow \qquad d_X(g \cdot x, g \cdot y) < \e .$$
  The   Cantor minimal system $(F_0,G_0,\Phi)$ associated to the fiber of a weak solenoid is equicontinuous.  

An \emph{automorphism} of $(X,G,\Phi)$ is a homeomorphism $h \colon X \to X$ which commutes with the   $G$-action on $X$.  That is, for every $x \in X$ and   $g \in G$,    $g \cdot h(x) = h(g \cdot x)$. We denote by ${Aut}(X,G,\Phi)$ the group of automorphisms of the action $(X,G,\Phi)$. Note that ${Aut}(X,G,\Phi)$ is a topological group using the compact-open topology on maps, and  is a closed subgroup of ${\it Homeo}(X)$. Given a homeomorphism  $\tau \colon  X \to Y$ which conjugates the actions  $(X,G,\Phi)$ and $(Y,G, \Psi)$, then $\tau$ induces a topological isomorphism $\tau^* \colon {Aut}(X,G,\Phi) \cong {Aut}(Y,G,\Psi)$. Thus, the properties of the group ${Aut}(X,G,\Phi)$ and its action on the space $X$ are topological conjugacy invariants of the Cantor minimal system.

An action $(X,G,\Phi)$ is \emph{homogeneous} if $Aut(X,G,\Phi)$ acts transitively on $X$, that is, for any $y \in X$ there is $h \in Aut(X,G,\Phi)$ such that $h(x) = y$. We note the following result of Auslander:

\begin{thm} \cite[Chapter 2, Theorem 13]{Auslander1988}\label{aus}
Let  $(X,G,\Phi)$ be a homogeneous Cantor minimal system. Then $(X,G,\Phi)$ is equicontinuous.
\end{thm}
The examples of Schori \cite{Schori1966}, Rogers and Tollefson in \cite{RT1971b}, and also those given later in this paper, show that the
statement of Theorem \ref{aus} does not have a converse:  there are equicontinuous Cantor minimal actions which are not homogeneous. These are the actions of interest for this work.

The standard technique for analyzing Cantor minimal systems with $G = \mZ$ is to introduce  nested sequences of Kakutani-Rokhlin partitions for the action, as used, for example, in the work of Hermann, Putnam and Skau \cite{HPS1992}. Forrest introduced in \cite{Forrest2000} an analogous construction of Kakutani-Rokhlin partitions  for $G = \mZ^n$ with $n > 1$. 
For group actions,   constructions of Kakutani-Rokhlin partitions were given in the work of  Gjerde and Johansen   \cite{GJ2000}, and  Cortez and Petite   \cite{CortezPetite2008,CortezPetite2014}.

 For an equicontinuous Cantor minimal system  $(X,G,\Phi)$,   there always exists a collection $\{\cP_i\}_{i \geq 0}$ of Kakutani-Rokhlin partitions,   as  discussed   in Section~\ref{sec-all-chains}. These partitions are constructed, for example, using the method of coding as described in  Appendix \ref{appendix}, following Clark and Hurder   \cite{ClarkHurder2013}.  The  elements of $G$ which fix   the initial clopen subset $V_i \subset X$ in the partition $\cP_i$ form a group $G_i$, called the \emph{isotropy group} of the action of $G$ at $V_i$. The assumption that  the partitions are nested implies  that  $G_{i+1} \subset G_i$, and thus $\{G_i\}_{i \geq 0}$ is a nested group chain,    as defined below.

\begin{defn}\label{gr-chain}
Let $G$ be a finitely generated group. A \emph{group chain} $\cG = \{G_i \mid i \geq 0\}$, with $G_0=G$, is a properly descending chain of subgroups of $G$, such that $|G:G_i| < \infty$ for every $i \geq 0$.
\end{defn}

Given a sequence of partitions $\{\cP_i\}_{i \geq 0}$, the intersection of the initial clopen subsets $V_i \subset X$ also defines a basepoint $x = \cap V_i$. Then  there exists a homeomorphism
  \begin{align*} 
  \phi \colon  X \to G_{\infty} = \lim_{\longleftarrow} \,\left\{G/G_i \to G/G_{i-1}\} = \{(G_0, g_1G_1, g_2G_2,\ldots) \mid g_jG_i = g_iG_i \textrm{ for all }j \geq i \right\},
  \end{align*}  
such that $\phi(x) = (eG_i)$, where $eG_i$ is the coset of the identity in $G/G_i$, and such that $\phi$ is equivariant with respect to the action of $G$ on $X$, where the natural left $G$-action on $G_{\infty}$ is  given by 
   \begin{align}\label{inv-limaction} \gamma \cdot (g_iG_i) = (\gamma g_i G_i), \textrm{ for all } \gamma \in G .\end{align}
    The dynamical system $(G_{\infty},G)$ is called the \emph{inverse limit dynamical system of a group chain} $\{G_i\}_{i \geq 0}$. Such an inverse limit dynamical system is completely determined by the group chain $\{G_i\}_{i \geq 0}$, and it has a canonical basepoint given by the sequence of cosets $(eG_i) \in G_{\infty}$.
    
The   group chain $\{G_i\}_{i \geq 0}$  associated to a minimal equicontinuous system $(X,G,\Phi)$ is not  unique, and in particular,     depends on the choice of $\{\cP_i\}_{i \geq 0}$ and the basepoint  point   $x \in X$. For this reason, we call $\{G_i\}_{i \geq 0}$ the \emph{group chain with basepoint} $x$ \emph{and partitions} $\{\cP_i\}_{i \geq 0}$. A fundamental point is to describe the dependence of the group chain on the partitions and on the basepoint chosen, and this relationship is made precise in Proposition~\ref{prop-uniqueness} in Section~\ref{sec-all-chains}, and uses the following two concepts.

For a given group $G$, let  $\fG$ denote the collection of all possible group chains in $G$. 
\begin{defn}\cite{FO2002,RT1971b}\label{defn-greq}
Let $G$ be a finitely generated group, and let $\{G_i\}_{i \geq 0}$ and $\{H_i\}_{i \geq 0}$ be   group chains in $\fG$, so that   $G_0=H_0=G$. Then 
\begin{enumerate}
\item The chains $\{G_i\}_{i \geq 0}$ and $\{H_i\}_{i \geq 0}$ are \emph{equivalent}, if and only if, there is a group chain $\{K_i\}_{i \geq 0}$ in $\fG$ and infinite subsequences $\{G_{i_k}\}_{k \geq 0}$ and $\{H_{j_k}\}_{k \geq 0}$ such that $K_{2k} = G_{i_k}$ and $K_{2k+1} = H_{j_k}$ for $k \geq 0$.
\item The chains $\{G_i\}_{i \geq 0}$ and $\{H_i\}_{i \geq 0}$ are \emph{conjugate equivalent}, if and only if, there exists a sequence $(g_i) \in G$ such that $g_iG_i = g_jG_i$ for all $i\geq 0$ and all $j \geq i$, so that  the group chains $\{g_iG_ig_i^{-1}\}_{i \geq 0}$ and $\{H_i\}_{i \geq 0}$ are equivalent. 
\end{enumerate}
\end{defn}

The dynamical meaning of the equivalences in Definition \ref{defn-greq} is given by the following theorem.

\begin{thm}\cite{FO2002}\label{equiv-rel-11}
Let $(G_{\infty}, G)$ and $(H_{\infty}, G)$ be inverse limit dynamical systems for group chains $\{G_i\}_{i \geq 0}$ and $\{H_i\}_{i \geq 0}$ in $\fG$. Then the following is true.
\begin{enumerate}
\item The group chains $\{G_i\}_{i \geq 0}$ and $\{H_i\}_{i \geq 0}$ are equivalent if and only if there exists a homeomorphism $\tau: G_{\infty} \to H_{\infty}$ equivariant with respect to the $G$-actions on $G_{\infty}$ and $H_{\infty}$, and such that $\phi(eG_i) = (eH_i)$.
\item The group chains $\{G_i\}_{i \geq 0}$ and $\{H_i\}_{i \geq 0}$ are conjugate equivalent if and only if there exists a homeomorphism $\tau: G_{\infty} \to H_{\infty}$ equivariant with respect to the $G$-actions on $G_{\infty}$ and $H_{\infty}$.
\end{enumerate}
\end{thm}

That is, an equivalence of two group chains corresponds to the existence of a \emph{basepoint-preserving} conjugacy between their inverse limit systems, while a conjugate equivalence of two group chains corresponds to the existence of a conjugacy between their inverse limit systems, which need not preserve the basepoint. Here is one application of Theorem~\ref{equiv-rel-11}.

\begin{cor}
Let $x \in X$ be a choice of basepoint. Then all group chains associated to a minimal system $(X,G,\Phi)$ with  basepoint $x$ are   equivalent. Moreover,   if $\{G_i\}_{i \geq 0}$ is a group chain for the action with basepoint $x$, and  $y \in X$, then the equivalence class   of group chains with   basepoint $y$ has a representative of the form $\{g_iG_ig_i^{-1}\}_{i \geq 0}$, where $\{g_i\} \in G$ with $g_iG_i = g_jG_i$ for all $j \geq 0$.
\end{cor}

The   two equivalence relations on group chains in Definition~\ref{defn-greq} were used  by Fokkink and Oversteegen in \cite{FO2002} to explain why the examples of Rogers and Tollefson in \cite{RT1971b} are   homogeneous. This requires the following notion, which is Definition~16 in \cite{FO2002}:

\begin{defn}\label{def-weaklynormal}
A  group chain $\{G_i\}_{i \geq 0}$ is \emph{weakly normal}, if   there exists some index $i_0 \geq 0$ such that the subchain $\{G_i\}_{i \geq i_0}$ defines a single   equivalence class. 
\end{defn}
This condition can be formulated in terms of the group chains as follows. The group chain $\{G_i\}_{i \geq 0}$ is weakly normal if there exists an index $i_0 \geq 0$ such that for the restricted action of $G_{i_0}$ on the invariant clopen set $V_{i_0}$, all group chains defined by this action for any basepoint $y \in V_{i_0}$ are equivalent.

Fokkink and Oversteegen proved the following result, which is Theorem~18 in \cite{FO2002}:
 \begin{thm}\label{thm-weaklynormal}
 Let $\{G_i\}_{i \geq 0}$ be the group chain associated to a weak solenoid $\cS$. 
Then $\{G_i\}_{i \geq 0}$ is a weakly normal  group chain if and only if $\cS$ is homogeneous.
 \end{thm}
It can be a difficult task to show that a given group chain $\{G_i\}_{i \geq 0}$ is  {weakly normal}, and the first point of this work is to obtain a computable invariant of a group chain which reveals whether the group chain is weakly normal, or not. Our approach is to recall a notion well-known in topological dynamics, and use this to define invariants of the conjugate  equivalence class of the group chain.

Associated to an equicontinuous Cantor minimal system  $(X,G,\Phi)$, there is a  compact profinite subgroup    $\oG$ which acts transitively on $X$. The group $\oG$ is the closure in the uniform topology of the image   subgroup  $\Phi(G)  \subset {\it Homeo}(X)$, and is a special case of the Ellis semigroup associated to the Cantor minimal system. This is discussed in detail in Section~\ref{sec-ellis}. For $x \in X$, let 
\begin{align}\label{iso-defn}\ds \oG_x = \{ \whg \in \oG \mid \whg(x) = x\}\end{align}
denote  the isotropy group at $x$ of the left  action of $\oG$ on $X$. 
Then  there is a natural identification $X \cong \oG/\oG_x$ of left $G$-spaces.  Note that for any $y \in X$ the isotropy groups $\oG_x$ and $\oG_y$ are conjugate by an element in $\oG$.

Proposition~\ref{prop-auto=normalizer} below shows that for the study of the orbits of ${Aut}(X,G,\Phi)$, it suffices to analyze  the normalizer subgroups in $\oG$ of the isotropy groups $\oG_x$ for $x \in X$. We  use the  \emph{calculus of  group chains} developed in Section~\ref{sec-all-chains}, to obtain   explicit representations for $\oG$ and $\oG_x$, which makes it possible to study the properties of these normalizer subgroups.
This is done in Section~\ref{sec-quotient} below, and the key results are summarized as follows.

Let $\{G_i\}_{i \geq 0}$ be a group chain  with partitions $\{\cP_i\}_{i \geq 0}$ and basepoint  point $x \in X$. For each $i \geq 0$, the \emph{core} of $G_i$ in $G$ is the normal subgroup ${\ds C_i = \bigcap_{g \in G} gG_ig^{-1}}$.
It is immediate  that    the collection $\{C_i\}_{i\geq 0}$ of cores forms a nested group chain of normal subgroups of $G$. Then the inverse limit 
\begin{equation}\label{eq-defcore}
C_{\infty} = \lim_{\longleftarrow} \,\left\{G/C_i \to G/C_{i-1} \right\}
\end{equation}
  is a profinite group, called the \emph{limit core}, and there is a natural coordinate action of $G$ on $C_{\infty}$, given by a formula analogous to \eqref{inv-limaction}. The \emph{discriminant group} $\cD_x$ of $(X,G,\Phi)$ at $x \in X$  is introduced in   Section~\ref{sec-quotient},  and is   a profinite subgroup of    $C_{\infty}$ such that   $C_{\infty}/\cD_x = G_{\infty}$. 
 Then we have:

 \begin{thm}\label{thm-quotientspace-1}
Let $(X,G,\Phi)$ be a minimal equicontinuous dynamical system, and let $\{G_i\}_{i \geq 0}$ be a group chain with basepoint $x \in X$  associated to the system.
Then there is a natural isomorphism of topological groups $ \oG \cong C_{\infty}$ such that $\oG_x \cong \cD_x$.
\end{thm}

Theorem~\ref{thm-quotientspace-1} can be interpreted as providing `coordinate representations' $C_{\infty}$ and $\cD_x$ for the   Ellis group $\oG$ of $(X,G,\Phi)$ and the isotropy groups $\oG_x$ for $x \in X$.  These `representations' depend on the  choice of partitions $\{\cP_i\}_{i \geq 0}$ with basepoint $x$, and are unique up to a topological isomorphism, as described in detail in Section \ref{sec-quotient}.

We now consider applications of Theorem~\ref{thm-quotientspace-1}. First, there is the following immediate corollary, which follows from the fact that for $x,y \in X$, the isotropy groups $\oG_x$ and $\oG_y$   are conjugate by an element of $C_{\infty}$. 

\begin{cor}
Let $(X,G,\Phi)$ be a minimal equicontinuous action, and let $\{G_i\}_{i \geq 0}$ be a group chain with partitions $\{\cP_i\}_{i \geq 0}$ and   basepoint $x$. Then the cardinality of $\cD_x$ is independent of the choice of the partitions and of the basepoint.
\end{cor}

Here is a basic application of  Theorem \ref{thm-quotientspace-1}, which is proved in Section \ref{sec-properties}.

\begin{thm}\label{cor-discrisnotnormal}
Let $(X,G,\Phi)$ be an equicontinuous Cantor minimal system,   $x \in X$ a basepoint, and $\oG_x$ the isotropy group of $x$. Then  
  $${\rm core}_{G}=  \bigcap_{k \in G} k \ \oG_x \ k^{-1}$$ 
 is the trivial group, and thus the maximal normal subgroup of $\oG_x$ is also the trivial group.
\end{thm}

It follows that $\oG_x$ is a normal subgroup of $\oG$ if and only if $\oG_x$ is trivial. Thus,  if   $\oG_x$ is non-trivial, then    $X \cong \oG/\oG_x$  does not have a   structure as a quotient group, and so $\oG_x$ can be viewed as a measure of non-homogeneity of the action on $X$. This leads to the following application of Theorems~\ref{thm-quotientspace-1} and \ref{cor-discrisnotnormal}:

\begin{cor}\label{cor-3-1}
Let $(X,G,\Phi)$ be an equicontinuous Cantor minimal system with a group chain $\{G_i\}_{i \geq 0}$ at $x$, and let $\cD_x$ be the discriminant group at $x$. Then $(X,G,\Phi)$ is homogeneous if and only if $\cD_x$ is trivial.
\end{cor}

Recall that a weak solenoid $\cS$ is regular if all of the coverings in a defining sequence for it are regular coverings. 
Similarly, a group chain $\{G_i\}_{i \geq 0}$ is \emph{regular} if each $G_i$ is normal in $G = G_0$.

Consider    the following    property of group chains:
 \begin{defn}\label{def-virtually-regular}
 Let $(X,G,\Phi)$ be a minimal equicontinuous action of a finitely generated group $G$ on a Cantor set $X$.  Then the action is \emph{virtually regular}  if there is a   subgroup $G_0' \subset G_0 = G$ of finite index   such that the restricted chain $\{G_i' = G_i \cap G_0'\}_{i \geq 0}$  is \emph{weakly normal} in $G_0'$. 
 \end{defn}
The normal  core of the subgroup $G_0' \subset G$ in Definition~\ref{def-virtually-regular} also satisfies the conditions of the definition, so we can assume without loss of generality that the subgroup $G_0'$  is normal in $G$.

 Note that a weakly normal chain is always virtually regular, but the converse need not be true. The   relation between these notions is discussed further in Remark~\ref{rmk-comparison}. Note that the   Rogers and Tollefson example in \cite{RT1972}, as given here in 
 Example~\ref{ex-RT}, is virtually regular but is not weakly normal.

Now suppose $\oG_x$ is non-trivial for some $x \in X$.
Since $\oG_x$ is a closed subgroup of a compact group $\oG$, there are two possibilities:   $\oG_x$ is finite, or $\oG_x$ is   a Cantor set.

In the case when $\cD_x$ is finite, Theorem~\ref{thm-core-finite} shows that the system $(X,G,\Phi)$ admits a special `nice group chain'. This is used to show the following    result, as shown in  Section~\ref{sec-weakly-reg-dx}.

 \begin{thm}\label{thm-virtually-regular}
 Let $(X,G,\Phi)$ be a minimal equicontinuous action of a finitely generated group $G$ on a Cantor set $X$. Suppose the discriminant group $\cD_x$ is finite for some $x \in X$, then $(X,G,\Phi)$  is virtually regular.
 \end{thm}

Theorem~\ref{thm-virtually-regular} has the following   application for weak solenoids, as discussed in Section \ref{sec-weakly-reg-dx}. It is also the motivation for the   work \cite{DHL2016b}.
 
 \begin{cor}\label{cor-almostcover}
 Let $M_0$ be a closed manifold, $x_0 \in M_0$ a fixed basepoint, and let $G = \pi_1(M_0 , x_0)$ be the fundamental group.
Let $f \colon \cS \to M_0$ be a weak solenoid, with  $X= f^{-1}(x_0)$ the Cantor fiber, and $x \in X$ a given basepoint. 
 Suppose the   monodromy action $(X,G,\Phi)$ on the fiber has   finite discriminant group $\cD_x$. Then there exists a finite-to-one covering map $p_0 \colon M_0' \to M_0$ such that the pullback solenoid $f' \colon  \cS' \to M_0'$ is homogeneous.
 \end{cor}
 
That is, for a weak solenoid whose holonomy action is    a virtually regular action, one obtains a homogeneous action by pulling back along a finite-to-one covering of the base manifold. Thus, the condition that $\cD_x$ is finite implies there is some finite covering of the weak solenoid which is homogeneous. A basic point is that this covering need not be related to one of the finite covers in the   inverse limit systems which define  $\cS$.  
 
Thus, if the  weak solenoid $\cS$ satisfying the hypotheses of Corollary~\ref{cor-almostcover} is not weakly normal, then it is not homogeneous by Theorem~\ref{thm-weaklynormal}. On the other hand, the pull-back solenoid $\cS'$ of Corollary~\ref{cor-almostcover} is homogeneous, hence cannot be homeomorphic to $\cS$.  These remarks yield the following consequence:
\begin{cor}\label{cor-distinctcoverngs}
 Let $M_0$ be a closed manifold with fundamental group $G$, and let $f \colon \cS \to M_0$ be a weak solenoid. Let $x_0 \in M_0$ be a basepoint, and $X= f^{-1}(x_0)$ be the Cantor fiber, with $x \in X$ the basepoint. Suppose the minimal equicontinuous action $(X,G,\Phi)$ on the fiber is not weakly normal, but has a non-trivial finite discriminant group $\cD_x$. Then there exists a finite-to-one covering map $p_0 \colon M_0' \to M_0$ such that the pull-back solenoid $\cS'$   is a covering of $\cS$   not homeomorphic to $\cS$. 
 \end{cor}

In the case when  the discriminant      $\cD_x$ is a Cantor group,   Corollary~\ref{cor-3-1} implies that the system $(X,G,\Phi)$ is not homogeneous.  Examples~\ref{main-6} and \ref{example63} construct  minimal Cantor actions which have infinite discriminant group.

The \emph{kernel} of a group chain $\{G_i\}_{i \geq 0}$    is the   subgroup
\begin{equation}
K_x = \bigcap_{i\geq 0}  ~ G_i ~ \subset ~ G \ .
\end{equation}
Here, the subscript on $K_x$ corresponds to the basepoint $x = (eG_i)$. When the group chain $\{G_i\}_{i \geq 0}$ is associated to a presentation of a weak solenoid $\cS$, with fiber $X =   f^{-1}(x_0)$,  then the kernel $K_x$ is identified with the fundamental group of the leaf $L_x \subset \cS$ containing the point $x$. The kernel thus has a natural geometric interpretation.

The kernel $K_x$ can either be a trivial group, a non-trivial normal subgroup of $G$, or a non-normal subgroup of $G$. Also, the kernel need not be finitely-generated, as in the case of the non-homogeneous solenoid constructed by  Schori in \cite{Schori1966}.
  Note that if  $\{H_i\}_{i \geq 0}$ is a group chain equivalent to $\ds \{G_i\}_{i \geq 0}$, then clearly $\ds \bigcap_i H_i = \bigcap_i G_i = K_x$.  

  On the other hand, given $y \in G_{\infty}$   let $\{H_i\}_{i \geq 0}$ be a group chain with basepoint $y$ which is conjugate equivalent to $\ds \{G_i\}_{i \geq 0}$, then its kernel $K_y$ need not equal $K_x$. 
  In  Section~\ref{sec-kernel} we give examples which show that   the kernel $K_x$ of a group chain $\{G_i\}_{i \geq 0}$ is not invariant under the conjugate equivalence of group chains.
  It seems to be a subtle property of a Cantor minimal system to determine whether,  given  a point $x \in X$ for which the kernel $K_x$ is  not a normal subgroup, hence is non-trivial, must the kernel $K_y$ be  non-trivial  for all conjugate chains.
  
  If $\ds \{G_i\}_{i \geq 0}$ is the group chain associated to a weak solenoid $\cS$, then 
the non-invariance of the kernel of a group chain under the change of basepoint can be formulated as a property of the holonomy groups  of leaves  of the foliation on $\cS$. The relationship between the non-homogeneity of a group chain defining a weak solenoid and the dynamics of $\cS$ are discussed further in the works  \cite{DHL2016b,DHL2016c}.

 Section~\ref{sec-properties} introduces the   notion of a group chain in \emph{normal form} with respect to the discriminant invariant $\cD_x$, as given by  Definition~\ref{def-normalform}, which can be viewed as weaker form of the conclusion of  Theorem~\ref{thm-core-finite}.  
Every  group chain is equivalent to one in  normal form  by Proposition~\ref{prop-specialchoice}. This result leads to the introduction of a new property of group chains, based on the relationship between the kernels of group chains and the discriminant group.

 \begin{defn}\label{def-stable}
Let $\{G_i\}_{i \geq 0}$ be a group chain in a normal form, and let $\cD_x$ be the discriminant group of the action $(G_{\infty},G)$ at $x = (eG_i)$. 
The  discriminant group $\cD_x$ is \emph{stable} if for any $y \in G_{\infty}$ every element in $\cD_y$ is represented by an element in the kernel $K_y$.
The group $\cD_x$ is \emph{unstable} if it is not stable.
\end{defn}

In Section~\ref{sec-kernel} the properties of kernels and their relationship to the discriminant group $\cD_x$ are considered.  
   In particular,    the discriminant group in Example~\ref{ex-RT} is infinite and unstable, while the discriminant group  in Example \ref{finitediscr} is finite and stable.

 The method of calculating the discriminant invariant for a group chain is   an effective technique for obtaining results about the dynamical properties of the equicontinuous systems thus obtained.  
Example~\ref{ex-RT}  calculates the discriminant group for the   Rogers and Tollefson example in \cite{RT1972}.

 Finally, in   Section~\ref{sec-examples}, we give   examples of minimal equicontinuous actions of non-abelian groups  with non-trivial discriminant group, and characterize whether they are weakly normal, virtually regular, and have unstable discriminant group.   These examples include the actions of the discrete Heisenberg group $\cH$, of the dihedral group, and of the generalized dihedral group.  
 
 \begin{remark}  {\rm
 The following table lists the examples and their properties. 
  
 \begin{center}
\begin{tabular}{cccccc}
 & type & weakly normal & virt. regular & \#~$\cD_x$ & stable \\
 \hline
 {\rm Example}~\ref{ex-RT}~ & dihedral & no & yes & $< \infty$  & no \\
  {\rm Example}~\ref{finitediscr}  ~ & product & yes & yes & $< \infty$ & yes \\
  {\rm Example}~\ref{eq-wrHeisenberg} ~ & Heisenberg & yes & yes & $< \infty$ & ? \\
  {\rm Example}~\ref{main-6}  ~ & Heisenberg & no & ? & $= \infty$ & no \\
  {\rm Example}~\ref{example63}  ~ & dihedral & yes & yes &  $= \infty$ & no \\
  {\rm Example}~\ref{example63-1} ~ & dihedral & yes & yes &  $= \infty$ & no
 \end{tabular} 
 \end{center}
}
\end{remark}
 
 \section{The Ellis group}\label{sec-ellis}

The concept  of the Ellis semigroup associated to a continuous group action $\Phi \colon G \times X \to X$   was introduced in the papers \cite{EllisGottschalk1960,Ellis1960}, and is treated in the books by Ellis \cite{Ellis1969,Ellis2014} and Auslander \cite{Auslander1988}.   The construction of $\whE(X,G,\Phi)$ is abstract, and it can be a   difficult problem to determine this group exactly, and the relation between  the algebraic properties of $\whE(X,G,\Phi)$ and dynamics of the action.
In this section, we recall some basic properties   of $\whE(X,G,\Phi)$, then consider the  results for the special case of equicontinuous Cantor minimal actions.

  Let $X$ be a compact Hausdorff  topological space, and let $X^X = {\it Maps}(X,X)$   have the topology of \emph{pointwise convergence on maps}. This is equivalent to the Tychonoff topology, where we identify  $X^X$ with the product space $\ds \prod_{x \in X} ~ X_x$  where $X_x = X$ for all $x \in X$. By the Tychonoff Theorem, $X^X$ is   a compact Hausdorff space.
Note that there is an inclusion ${\it Homeo}(X) \subset X^X$, so   each $g \in G$ defines an element   $\whg \in X^X = {\it Maps}(X,X)$. Let
\begin{equation}
\whG = \{ \whg  \mid g \in G \} \subset {\it Homeo}(X) \subset X^X
\end{equation}
which is a group under the natural composition of elements of $X^X$.  The action $\Phi$ defines a map, $\Phi \colon G \to \whG$ with $\Phi(g) = \whg$, whose kernel is a normal subgroup, denoted by $L_x \subset G$. 
\begin{defn}
The \emph{Ellis (enveloping) semigroup} of    a continuous group action $\Phi \colon G \times X \to X$  is the closure $\whE(X,G,\Phi) = \overline{\whG} \subset X^X$ in the topology of \emph{pointwise convergence} on maps.
\end{defn}
Ellis proved that $\whE(X,G,\Phi)$ forms a right topological semigroup.  
As $\whE(X,G,\Phi)$ is a closed subset of a compact Hausdorff space, it is a compact set. 
However,   $\whE(X,G,\Phi)$ need not be separable. 
 If all the elements of $\oG$ are defined as the sequential limit of elements of $\whG$, then the dynamical system  is said to be \emph{tame},  in the sense of Glasner \cite{Glasner2006}. See Glasner \cite{Glasner2007b} for a detailed description of tame systems. Otherwise, the elements of $\oG$ must be  defined using ultrafilters, and it can be extremely difficult to give a detailed description of this group \cite{Glasner2007a}.

In the case of an equicontinuous  Cantor minimal system $(X,G,\Phi)$, the properties of its   Ellis semigroup $\whE(X,G,\Phi)$ are   discussed by Ellis in \cite{Ellis1960}, Glasner   in \cite[Proposition~2.5]{Glasner2007a},
and by Auslander in \cite[Chapter 3]{Auslander1988}. The equicontinuity assumption on the action implies that the pointwise and uniform topologies on $\whE(X,G,\Phi)$ coincide, and thus each element of $\whE(X,G,\Phi)$ is the sequential limit of points in $\whG$, so that $(X,G,\Phi)$ is a tame dynamical system. Consequently, the space $\whE(X,G,\Phi)$ is separable, and moreover we have the basic result:

\begin{thm} \cite[Chapter 3]{Auslander1988}.
Let $(X,G,\Phi)$ be an equicontinuous Cantor minimal system. Then $\whE(X,G,\Phi)$ is a compact topological group, which is naturally   identified with the closure  $\oG$  of $\Phi(G)  \subset {\it Homeo}(X)$ in the \emph{uniform topology on maps}.
\end{thm}
   
   In the following, we identify $\oG = \whE(X,G,\Phi)$. As each element of $\oG$ is the limit of a sequence of points in $\whG$, we use the notation $(g_i)$ to signify a sequence $\{g_i \mid i \geq 1\} \subset G$ such that the sequence $\{\whg_i = \Phi(g_i) \mid i \geq 1\} \subset {\it Homeo}(X)$ converges in the uniform topology. 
   
Let   $\whPhi \colon \oG \times X \to X$ denote the action of $\oG$ on $X$.
   Note that for $(g_i) \in \oG$, the action   $\whPhi(g_i) \colon X \to X$, for $x \in X$, is given by  $\whPhi(g_i)(x) = (g_i) \cdot x = \lim \ g_i \cdot x= \lim \ \Phi(g_i)(x)$.

   The orbit of any point $x \in X$ under the action of $\oG$ is a compact subset of $X$. By the minimality  assumption,  the orbit $\whG(x)$ is dense in $X$ for any $x \in X$, so it follows that the orbit $\oG(x) = X$ for any $x \in X$. That is, the group $\oG$ acts transitively on $X$. Then for the isotropy group of $x$,  
\begin{align}\label{iso-defn2}
 \oG_x = \{ (g_i) \in \oG \mid (g_i) \cdot x = x\},
\end{align}
  we have the natural identification $X \cong \oG/\oG_x$ of left $G$-spaces. 
The pair $(\oG, \oG_x)$ is said to ``classify'' the   Cantor minimal system $(X,G,\Phi)$, in the sense of topological dynamics.

 Since the action of $\oG$ is transitive on $X$, given $y \in X$, then there exists $(g_i) \in \oG$ such that $(g_i) \cdot x  =y$. It follows that 
  \begin{align}\label{eq-profiniteconj}
  \oG_y = (g_i) \cdot \oG_x \cdot (g_i)^{-1}.  
  \end{align}
This identity implies the following fact.

\begin{lemma}\label{cardinality-conjuagte}
Let $(X,G,\Phi)$ be an equicontinuous Cantor minimal system. Then
the cardinality of $\oG_x$ is independent of the point $x \in X$; that is, for any $y \in X$ we have
  \begin{align}\label{card-statement} card (\oG_y) = card (\oG_x). \end{align}
\end{lemma}

We also have the basic observation about closed subgroups of a Cantor group:
\begin{lemma}\label{cardinality-lemma}
Let $(X,G,\Phi)$ be an equicontinuous Cantor minimal system, and let $x \in X$. Then
the isotropy subgroup $\oG_x$  is either a finite group, or is a Cantor group.
\end{lemma}
\proof The groups   $\oG_x$ and $\oG$ are compact, so if $\oG_x$ is not finite, then it contains a limit point of itself. As $\oG_x$ is a group, every point must then be a limit point. Since $\oG$ is totally disconnected, the same holds for $\oG_x$ and thus it is also a Cantor group.
\endproof

 Next, we consider the relations between the above constructions and the group of automorphisms ${Aut}(X,G,\Phi)$ of an equicontinuous  Cantor minimal system $(X,G,\Phi)$.
First, we show:
\begin{lemma}\label{lem-auto}
${Aut}(X,G,\Phi) = {Aut}(X,\oG, \whPhi)$.  
\end{lemma}
\proof
The inclusion ${Aut}(X,G,\Phi) \supseteq {Aut}(X,\oG, \whPhi)$ follows from the inclusion $G \subset \oG$.

To show the reverse inclusion, let  $(g_i) \in \oG$.  
Let $h \in {Aut}(X,G,\Phi)$ and $y \in X$. The map $h$ is a homeomorphism of the compact set $X$, so  we have
\begin{equation}\label{eq-approxmap}
h((g_i) \cdot y) =   h(\lim \, g_i \cdot y) = \lim \, g_i \cdot h(y) = (g_i) \cdot h(y).
\end{equation}
As \eqref{eq-approxmap} holds for all $y \in X$, we conclude that $h \in {Aut}(X,\oG, \whPhi)$, as was to be shown.
\endproof

Now let  $\whcN(\oG_x) \subseteq \oG$ denote the normalizer of $\oG_x$ in $\oG$. Then the quotient group    
\begin{equation}\label{eq-nolmalizedauto}
\whcA(X,\oG,x) \equiv \whcN(\oG_x)/\oG_x
\end{equation}
acts   \emph{on the right} on the coset space $\oG/\oG_x$. The right action of $\whcA(X,\oG,x)$ commutes with the left action of $\oG$ on $\oG/\oG_x$, so there is an   inclusion  $\whcA(X,\oG,x) \subseteq {Aut}(X,\oG,\whPhi)$.
The following result is the analog  of a well-known result for finite  covering spaces   \cite[Corollary~7.3]{Massey1967}.

\begin{prop}\label{prop-auto=normalizer}
   ${Aut}(X,G,\Phi)$ is isomorphic to    $\whcA(X,\oG,x)$. 
\end{prop}
\proof
By the above remarks, it suffices to show ${Aut}(X,G,\Phi) \subseteq \whcA(X,\oG,x)$. 

Let $h \in {Aut}(X,G,\Phi)$. By Lemma~\ref{lem-auto}, the homeomorphism $h \colon X \to X$ also commutes with the action of  $\oG$ on $X$.  Let $z = h(x)$. Then for $(g_i) \in \oG_x$ we have $(g_i) \cdot x = x$,   so by the uniform convergence of the maps $\{\Phi(g_i)\}$, we have
$$(g_i)  \cdot z = (g_i)  \cdot h(x) = h((g_i)  \cdot x) =   h(x)   = z ,$$
so that $\oG_x \subseteq \oG_z$. The reverse inclusion holds by considering the map $h^{-1}$, so we obtain 
$\ds \oG_x = \oG_z$. 

  Now let $(k_i) \in \oG$ such that $(k_i) \cdot x = z$. Then by \eqref{eq-profiniteconj} we have 
\begin{equation}\label{eq-normalized}
\oG_x = \oG_z = (k_i) \cdot \oG_x \cdot (k_i)^{-1} 
\end{equation}
which implies that $(k_i) \in \whcN(\oG_x)$, and hence defines a right action on $X$. Denote this right action by 
$\rho(k_i) \in {\it Homeo}(X)$ where $\rho(k_i)(x) = x \cdot (k_i)$.

As $x$ is the coset of the identity, that is $x = e \cdot \oG_x$ under the identification $X \cong \oG/\oG_x$ of left $G$-spaces, we obtain that  the right action of $(k_i)$ satisfies $\rho(k_i)(x) = x \cdot (k_i) = z = h(x)$.

 Then consider the composition $\psi =    \rho(k_i) \circ h^{-1} \in {\it Homeo}(X)$ which  fixes the point $x$. By the remarks above, we have that  $\psi \in {Aut}(X,G,\Phi)$. Then for any $g \in G$, set $y = g \cdot x$, then we have 
  $$ y =g \cdot x = g \cdot \psi(x) = \psi(g \cdot x) = \psi(y),$$ 
  so that $\psi$ also fixes every point in the $G$-orbit of $y$. As the orbit $G \cdot x$ is dense in $X$, $\psi$ must be the identity. That is,  $h = \rho(k_i)  \in \whcA(X,\oG,x)$, as was to be shown.
\endproof

We obtain from Proposition~\ref{prop-auto=normalizer}   the following property of  ${Aut}(X,G,\Phi)$:
 
\begin{cor}\label{cor-rephr-2}
Let $(X,G,\Phi)$  be an equicontinuous minimal action of a finitely generated group $G$ on a Cantor set $X$. Then then number of orbits of ${Aut}(X,G,\Phi)$ is equal to the index of $\whcN(\oG_x)$  in $\oG$.
\end{cor}

For the actions where the automorphism group has a finite number of orbits, we have the following.

\begin{prop}\label{WR-clopen}
Let $(X,G,\Phi)$  be an equicontinuous minimal action of a finitely generated group $G$ on a Cantor set $X$, and suppose $|\oG:\whcN(\oG_x)| < \infty$. Then the orbits of ${Aut}(X,G,\Phi)$ are closed and open subsets of $X$.
\end{prop}

\proof  As $\whcN(\oG_x)$ is a closed subgroup of $\oG$, hence the orbits of its right action on $X$ are closed. Since ${Aut}(X,G,\Phi)$ has a finite number of orbits, then every orbit is an open subset of $X$, being the complement of a finite union of closed sets. Thus the orbits of ${Aut}(X,G,\Phi)$ are closed and open subsets of $X$. \endproof

For an equicontinuous Cantor minimal system, Proposition~\ref{prop-auto=normalizer} shows that the study of the group ${Aut}(X,G,\Phi)$, and in particular whether it acts transitively on $X$, is equivalent to the study of the normalizer group $\whcN(\oG_x)$ for the choice of a basepoint $x \in X$. In the following sections, 
 we will develop effective tools for the calculation of the groups $\oG$, $\oG_x$ and $\whcN(\oG_x)$ from a group chain associated with the action.

\section{A `calculus' of group chains}\label{sec-all-chains}

We recall some preliminary results about the group chains associated to   equicontinuous Cantor minimal actions. The main   reference for these results is the paper by  Fokkink and Oversteegen \cite{FO2002}, and Clark and Hurder \cite{ClarkHurder2013} for the construction of partitions associated to an action.

The following is the basic result which relates the action on $X$ with group chains in $G$.

\begin{prop}\label{prop-AFpres}
Let $(X,G,\Phi)$ be an   equicontinuous Cantor minimal system, and let $x \in X$ be a point. Then there exists an infinite sequence of closed and open sets $X = V_0 \supset V_1\supset V_2\supset \cdots$ with $  \{x\} = \bigcap_{i} V_i$, which have the following properties.
\begin{enumerate}
\item \label{property-1} For each $i \geq 1$ the collection $\cP_i =\{g \cdot V_i\}_{g \in G}$ is a finite partition of $X$ into clopen sets. 
\item \label{property-2} We have ${ diam}(g \cdot V_i ) < 2^{-i}$ for all $g \in G$ and all $i \geq 0$.
\item The collection of elements which fix $V_i$, that is,
 $$G_i = \{g \in G \mid g \cdot V_i= V_i\},$$
 is a subgroup of finite index in $G$. More precisely, $|G:G_i| = { card}(\cP_i)$.
\item \label{property-5} There is a homeomorphism $\ds \phi \colon  X \to G_{\infty} = \lim_{\longleftarrow} \,\left\{G/G_{i+1} \to G/G_i \right\}$ equivariant with respect to the action of $G$ on $X$ and $G_{\infty}$, which maps $x$ onto $(eG_i) \in G_{\infty}$.
\end{enumerate}
\end{prop}

We note that property \eqref{property-1} implies that the sets in $\cP_i$ are permuted by the action of $G$, that is, if $U \in \cP_i$ and $g \in G$, then $g \cdot U = V$ for some $V \in \cP_i$.

The proof of Proposition \ref{prop-AFpres} uses the method of \emph{coding orbits}, described in detail in \cite[Section 6]{ClarkHurder2013}, where it was used to obtain a similar statement for a more general setting of equicontinuous \emph{pseudo}group actions. For completeness, we present the proof of Proposition \ref{prop-AFpres} in Appendix \ref{appendix}.

In Proposition \ref{prop-AFpres}, note that since $V_{i+1} \subset V_i$, then there is an inclusion of isotropy subgroups $G_{i+1} \subset G_i$, and so the collection $\{G_i\}_{i \geq 0}$ is a nested chain of subgroups of finite index. From Proposition \ref{prop-AFpres}, it is clear that the group chain $\{G_i\}_{i \geq 0}$ depends on the choice of partitions $\{\cP_i\}_{i \geq 0}$ and on the choice of $x \in X$, so we call $\{G_i\}_{i \geq 0}$ the \emph{group chain at }$x$ \emph{with partitions} $\{\cP_i\}_{i \geq 0}$.

The dependence of a group chain on the choice of the partitions and of the basepoint is described in the following proposition. Recall that Definition~\ref{defn-greq} introduced the notion of equivalent chains.
 
\begin{prop}\label{prop-uniqueness}
Let $(X,G,\Phi)$ be a minimal equicontinuous action, and let $\{G_i\}_{i \geq 0}$ be a group chain with partitions $\{\cP_i\}_{i \geq 0}$ and basepoint $x \in X$. Then the following hold:
\begin{enumerate}
\item Let $\{H_i\}_{i \geq 0}$ be a group chain with the same basepoint $x$ but different collection of partitions $\{\cQ_i\}_{i \geq 0}$. Then $\{H_i\}_{i \geq 0}$ is equivalent to $\{G_i\}_{i \geq 0}$.
\item Let $y \in X$, and let $\{H_i\}_{i \geq 0}$ be a group chain at $y$ with partitions $\{\cP_i\}_{i \geq 0}$. Then for every $i \geq 0$, there is an element $g_i \in G$ such that 
   $$\{H_i\}_{i \geq 0} = \{g_iG_ig_i^{-1}\}_{i \geq 0}.$$
 In particular, $g_i G_i = g_jG_i$ for all $j \geq 0$.  
\item \label{part3-uniqueness} Let $y \in X$ and let $\{H_i\}_{i \geq 0}$ be a group chain at $y$ with partitions $\{\cQ_i\}_{i \geq 0}$. Then there exists a sequence of elements $(g_i)_{i\geq 0}$ in $G$ such that
$\{H_i\}_{i \geq 0}$  is equivalent to $\{g_iG_ig_i^{-1}\}_{i \geq 0}$. That is, $\{G_i\}_{i \geq 0}$   is conjugate equivalent to $\{H_i\}_{i \geq 0}$.
\end{enumerate}
\end{prop}

\proof The result of Proposition \ref{prop-uniqueness} seems to be well-known; we give a proof for completeness.

  $(1).$ By Proposition \ref{prop-AFpres}, the group chain $\{G_i\}_{i \geq 0}$ is given by a collection of finite partitions $\{\cP_i\}_{i \geq 0}$, such that 
 $\cP_i = \{g \cdot V_i\}_{g \in G}$, where $V_i$ is a clopen set in $X$, such that $V_{i+1} \subset V_i$, $G_i$ is the isotropy group of $V_i$, we have $\bigcap V_i = \{x\}$, and $diam (g \cdot V_i) < 2^{-i}$ for $i \geq 0$.

Similarly, the group chain $\{H_i\}_{i \geq 0}$ is given by a collection of finite partitions $\{\cQ_i\}_{i \geq 0}$, such that $\cQ_i = \{g \cdot U_i \mid g \in G\}$, where $U_i$ is a clopen set in $X$, such that $U_{i+1} \subset U_i$, the group $H_i$ is the isotropy group of $U_i$, we have $\bigcap U_i = \{x\}$, and $diam (g \cdot U_i) < 2^{-i}$ for $i \geq 0$.

Since the diameters of the sets in $\{\cP_i\}$ and $\{\cQ_i\}$ tend to $0$ as $i \to \infty$, for every $i \geq 0$ we can find $j \geq i$ such that $U_j \subseteq V_i$ and $V_j \subseteq U_i$. 

Then the inclusion $U_j \subseteq V_i$ implies that $H_j \subseteq G_i$, and the inclusion $V_j \subseteq U_i$ implies that $G_j \subseteq H_i$ for the isotropy subgroups of the sets. Therefore, one can choose subsequences of clopen sets $\{V_{i_k}\}_{k \geq 0}$ and $\{U_{j_k}\}_{k \geq 0}$ satisfying $V_{i_1} \supset U_{j_2} \supset V_{i_k} \supset \ldots$, which implies that there is a nested group chain
\begin{equation}\label{eq-groupequiv}
G_{i_1} \supset H_{j_1} \supset G_{i_2} \supset H_{j_2} \supset \ldots,
\end{equation}
and so $\{H_i\}_{i \geq 0}$ and $\{G_i\}_{i \geq 0}$ are equivalent.

$(2).$  Now let $\{H_i\}_{i \geq 0}$ denote the group chain at $y \in X$ with partitions $\{\cP_i\}_{i \geq 0}$. Note that $y \in g_i \cdot V_i$ for some $g_i \in G$, by property \eqref{property-1} in Proposition \ref{prop-AFpres}. Then the isotropy group of the action $(X,G,\Phi)$ at $y$ is $H_i = g_i G_ig_i^{-1}$. Repeating this procedure for all $i \geq 0$, we obtain that $\{H_i\}_{i \geq 0} = \{g_i G_i g_i^{-1}\}_{i \geq 0}$. The collection $\{g_i\} \in G$ satisfies $g_iG_i = g_jG_i$ for all $j \geq i$, since $\{H_i\}_{i \geq 0}$ is nested. 

The proof of $(3)$ is just a combination of the proofs for $(1)$ and $(2)$.
\endproof

The dynamical meaning of group chain equivalences is then described by the following theorem.

\begin{thm}\label{thm-equiv-rel}
Let $(G_{\infty}, G)$ and $(H_{\infty}, G)$ be inverse limit dynamical systems for group chains $\{G_i\}_{i \geq 0}$ and $\{H_i\}_{i \geq 0}$. Then we have:
\begin{enumerate}
\item  The group chains $\{G_i\}_{i \geq 0}$ and $\{H_i\}_{i \geq 0}$ are equivalent if and only if there exists a homeomorphism $\tau: G_{\infty} \to H_{\infty}$ equivariant with respect to the $G$-actions on $G_{\infty}$ and $H_{\infty}$, and such that $\tau(eG_i) = (eH_i)$.
\item The group chains $\{G_i\}_{i \geq 0}$ and $\{H_i\}_{i \geq 0}$ are conjugate equivalent if and only if there exists a homeomorphism $\tau: G_{\infty} \to H_{\infty}$ equivariant with respect to the $G$-actions on $G_{\infty}$ and $H_{\infty}$.
\end{enumerate}
\end{thm}

That is, equivalence of group chains corresponds to the existence of a \emph{basepoint-preserving} homeomorphism between their respective inverse limit dynamical systems. If two chains are conjugate equivalent, then their inverse limit dynamical systems are conjugate, but the conjugating homeomorphism need not preserve the basepoint.

Theorem~\ref{thm-equiv-rel} is proved in \cite{FO2002} in the setting of weak solenoids using path lifts. We sketch here an alternative proof by means of standard group theory.

\proof $(1)$. First note that a group chain $\{G_i\}_{i \geq 0}$ is trivially equivalent to its subsequence $\{G_{i_k}\}_{k \geq 0}$, and it is immediate that the corresponding inverse limit systems $(G_{\infty},G, (eG_i))$ and $(G_{\infty,s},G,(eG_{i_k}))$ are conjugate by a basepoint preserving homeomorphism.

If $\{G_i\}_{i \geq 0}$ and $\{H_i\}_{i \geq 0}$ are equivalent group chains, then there are subsequences $\{G_{i_k}\}_{k \geq 0}$ and $\{H_{j_k}\}_{k \geq 0}$, and a group chain $\{K_i\}_{i \geq 0}$, such that $K_{2k} = G_{i_k}$ and $K_{2k+1} = H_{j_k}$ for $k \geq 0$. Then the conjugacy of $G_{\infty}$ and $H_{\infty}$ by a basepoint preserving homeomorphism follows from the existence of the conjugacies between the inverse limit dynamical systems corresponding to the pairs $\{G_i\}_{i \geq 0}$ and $\{G_{k_i}\}_{i \geq 0}$, $\{G_{k_i}\}_{i \geq 0}$ and $\{K_i\}_{i \geq 0}$, $\{K_i\}_{i \geq 0}$ and $\{H_{k_i}\}_{i \geq 0}$, and $\{H_{k_i}\}_{i \geq 0}$ and $\{H_i\}_{i \geq 0}$.

For the converse, suppose there exists a homeomorphism $\tau \colon   G_{\infty} \to H_{\infty}$, equivariant with respect to the $G$-action, and such that $\tau(eG_i) = (eH_i)$. 

Let $\delta_i \colon   H_{\infty} \to G/H_i$, and let $\tilde{\tau}_i = \delta_i \circ \tau \colon   G_{\infty} \to G/H_i$. Also, let  $\theta_i \colon   G_{\infty} \to G/G_i$ and $\tilde{\nu}_i = \theta_i \circ \tau^{-1} \colon  H_{\infty} \to G/G_i$. Note that since $\tau(eG_i) = (eH_i)$, we have $\tilde{\tau}_i(eG_i) = eH_i$, and $\tilde{\nu}_i(eH_i) = eG_i$.

Since $G/H_i$ is a finite space, by \cite[Lemma 1.1.16]{RZ2000} there exists $k_i>0$ and a surjective map $\tau_{i} \colon  G/G_{k_i} \to G/H_i$ such that $\tilde{\tau}_i = \tau_i \circ \theta_{k_i}$. Since $\theta_{k_i}(eG_i) = eG_{k_i}$, and $\tilde{\tau}_i(eG_i) = eG_i$, then $\tau_i(eG_{k_i}) = eH_i$, and so $G_{k_i} \subset H_i$.
By a similar argument, applied to $G/G_{k_i}$, there exists $\ell_i > 0$ such that $H_{\ell_i} \subset G_{k_i}$. Inductively, we obtain a sequence 
  $$ H_1 \supset G_{k_1} \supset H_{\ell_1} \supset G_{k_{\ell_1}} \supset \cdots,$$
  which shows that $\{G_i\}_{i \geq 0}$ and $\{H_i\}_{i \geq 0}$ are  equivalent.
  
 $(2)$. The following lemma says that if we change the basepoint, then the resulting system is conjugate to the given one. 
\begin{lemma}\label{thm-conjchains}
Let $\{G_i\}_{i \geq 0}$ be a group chain, and let $\{H_i\}_{i \geq 0} = \{g_iG_ig^{-1}_i\}_{i \geq 0}$, where $g_jG_i = g_i G_i$ for any $j \geq i$.  Let 
\begin{eqnarray*}
 G_{\infty} & = &   \lim_{\longleftarrow}\, \left\{\theta^i_{i-1} \colon  G/G_i \to G/G_{i-1}\right\},\\
 H_{\infty} & = & \lim_{\longleftarrow}\, \left\{\delta^{i}_{i-1} \colon  G/H_i \to G/H_{i-1}\right\}
\end{eqnarray*}
be the inverse limits spaces. Then there exists a homeomorphism $\tau \colon   G_{\infty} \to H_{\infty}$ such that $\tau(eG_i) = (g_i^{-1}H_i)$, and such that $\tau(s \cdot (h_i)) = s \cdot (\tau(h_i))$ for all $(h_i) \in G_{\infty}$ and all $s \in G$.
\end{lemma}

\proof Let $X = G_{\infty}$, and $x = (eG_i)$. For each $i \geq 0$, define 
  $$V_i = \{(g_iG_i) \in G_{\infty} \mid g_kG_k = eG_k \textrm{ for }k \leq i \}.$$
 Then 
  $$g \cdot V_i = \{ (h_iG_i) \in G_{\infty} \mid h_kG_k = g G_k \textrm{ for }k \leq i\},$$
and $\cP_i = \{g \cdot V_i\}_{g \in G}$, $i \geq 0$, is a collection of partitions corresponding to the group chain $\{G_i\}_{i \geq 0}$ with basepoint $x = (eG_i)$, and there is a bijective map $\tilde{\phi}_i \colon   \cP_i \to G/G_i$ defined by $\tilde{\phi}_i(g \cdot V_i) = gG_i$.

Since the chain $\{H_i\}_{i \geq 0} = \{g_iG_ig^{-1}_i\}_{i \geq 0}$ is nested, the diameter of $g_i \cdot V_i$ tends to $0$, and $G_{\infty}$ is compact, the intersection $\ds  \{y\} = \bigcap_{i} \, g_i \cdot V_i$ is non-empty. 

For each $i \geq 0$,  set  $ V_i ' = g_i \cdot V_i$, which we consider as the ``base partition'' and define the corresponding partition  $\cP_i' = \{ h \cdot V_i'\}_{h \in G}$ of $X$.   Then the partitions $\cP_i'$ and $\cP_i$ differ only by the ordering of sets,  that is, there is a bijection $\jmath_{\alpha_i} \colon   \cP_i \to \cP_i'$, induced by a permutation $\alpha_i \colon   (1,2, \ldots, \kappa_i) \to (1,2,\ldots,\kappa_i)$, such that $V_i^k = V_i^{\alpha_i(k)}$ as subsets of $G_{\infty}$.

The isotropy group of $g_i \cdot V_i = V_i'$ is $H_i = g_iG_ig_i^{-1}$, and there is a bijective map $\tilde{\psi_i} \colon   \cP_i' \to G/H_i$ for all $i \geq 0$. Define a bijective map $\tau_i \colon  G/G_i \to G/H_i$ by 
  $$\tau_i = \tilde{\psi}_i \circ \jmath_{\alpha_i} \circ (\tilde{\phi}_i)^{-1},$$ 
that is, $hG_i$ and $\tau_i(hG_i)$ correspond to the same subset of $G_{\infty}$. In particular, $\tau_i(eG_i) = g_i^{-1}H_i$. 

Since the action of $G$ permutes the sets in the partition $\cP_i$, the map $\tau_i$ is equivariant with respect to the $G$-actions on $G/G_i$ and on $G/H_i$. The mappings $\tau_i$ are compatible with the bonding maps of $G_{\infty}$ and $H_{\infty}$, and so there is an induced map $\tau \colon   G_{\infty} \to H_{\infty}$ which is clearly a bijection and is equivariant with respect to the action of $G$ on $G_{\infty}$ and $H_{\infty}$. By definition, $\tau(eG_i) = (g_i^{-1}H_i)$.
\endproof

We next give the proof of part (2) of Theorem~\ref{thm-equiv-rel}. 
Suppose $\{G_i\}_{i \geq 0}$ and $\{H_i\}_{i \geq 0}$ are conjugate equivalent, then $\{H_i\}_{i \geq 0}$ is equivalent to $\{g_i G_i g_i^{-1}\}_{i \geq 0}$ for some sequence $(g_i)_{i\geq 0}$ in $G$ such that $g_iG_i = g_jG_i$ for all $j \geq i$.
Denote by $(G_{\infty}',G)$ the inverse limit dynamical system of the group chain $\{g_i G_i g_i^{-1}\}_{i \geq 0}$. By Lemma \ref{thm-conjchains} the systems $(G_{\infty}',G)$ and $(G_{\infty},G)$ are conjugate, and by $(1)$ the systems $(G_{\infty}',G)$ and $(H_{\infty},G)$ are conjugate by a basepoint-preserving homeomorphism. Then $(G_{\infty},G)$ and $(H_{\infty},G)$ are conjugate.

For the converse, suppose $\tau \colon   G_{\infty} \to H_{\infty}$ is a homeomorphism equivariant with respect to the $G$-actions on $G_{\infty}$ and $H_{\infty}$. Let $x = (eG_i)$, and let $y = (g_iG_i) = \tau^{-1}(eH_i)$. 

By a procedure similar to the one in Lemma \ref{thm-conjchains}, given $\{G_i\}_{i \geq 0}$, we can define the partitions $\{\cP_i\}_{i \geq 0}$  of $G_{\infty}$, and given $\{H_i\}_{i \geq 0}$, we can define the partitions $\{\cQ_i\}_{i \geq 0}$ of $H_{\infty}$. 

Then $\cQ_i' = \tau^{-1}(\cQ_i)$ is a sequence of partitions of $G_{\infty}$. Since $\tau$ is equivariant with respect to the $G$-actions, the sets of $\cQ_i'$ are permuted by the action of $G$ on $G_{\infty}$, and the isotropy group of a set $U_i \in \cQ_i'$, containing $y$, is $H_i$. Then $\{H_i\}_{i \geq 0}$ is a group chain at $y$ with partitions $\{\cQ_i'\}_{i \geq 0}$. Then by Proposition \ref{prop-uniqueness} the group chains $\{G_i\}_{i \geq 0}$ and $\{H_i\}_{i \geq 0}$ are conjugate equivalent.
\endproof

We note that a basepoint preserving homeomorphism of $(X,G,\Phi)$ onto itself is precisely an automorphism of $(X,G,\Phi)$, and so 
Theorem~\ref{thm-equiv-rel} has the following corollary.

\begin{cor}\label{equiv-to-auto}
Let $(X,G,\Phi)$ be an equicontinuous minimal system, and let $\{G_i\}_{i \geq 0}$ be a group chain with partitions $\{\cP_i\}_{i \geq 0}$ at $x \in X$, and let $\{H_i\}_{i \geq 0}$ be a group chain with partitions $\{\cQ_i\}_{i \geq 0}$ at $y \in X$. Then there exists $h \in {Aut}(X,G,\Phi)$ such that $h(x) = y$ if and only if $\{G_i\}_{i \geq 0}$ and $\{H_i\}_{i \geq 0}$ are equivalent group chains.
 \end{cor}
 
Recall that Definition~\ref{defn-greq}  defined  equivalence  and conjugate equivalence  relations on  the collection $\fG$ of all group chains in $G$.    For an action $(X,G,\Phi)$, let $\{G_i\}_{i \geq 0}$ be a group chain with partitions $\{\cP_i\}_{i \geq 0}$ and basepoint $x$. 
Let  $\fG(\Phi)$ denote the class of group chains in $G$ which are \emph{conjugate equivalent} to $\{G_i\}_{i \geq 0}$.
    Theorem~\ref{thm-equiv-rel} yields the following:

\begin{cor}\label{cor-conjequiv}
For an action $(X,G,\Phi)$, let $\{G_i\}_{i \geq 0}$ be a group chain with partitions $\{\cP_i\}_{i \geq 0}$ and basepoint $x$.
 Then a group chain $\{H_i\}_{i \geq 0}$  is in $\fG(\Phi)$ if and only if there exists a collection of partitions $\cQ_i = \{g \cdot U_i\}_{g \in G}$, where $U_i$ is a clopen set, and ${\bigcap_i U_i = \{y\} \subset X}$, such that $H_i$ is the isotropy group of the action at $U_i$, for all $i \geq 0$.
\end{cor}

That is, the class $\fG(\Phi)$ of group chains   which are conjugate equivalent to  a given group chain $\{G_i\}_{i \geq 0}$,  contains all possible chains which can be produced for the action $(X,G,\Phi) = (G_{\infty},G)$ by the choice of changing partitions $\{\cP_i\}_{i \geq 0}$ and   basepoint $x \in X$.
Then we have:
\begin{thm}\cite{FO2002}\label{thm-rephr}
An equicontinuous Cantor minimal system  $(X,G,\Phi)$   is homogeneous  if and only if all group chains in $\fG(\Phi)$ are equivalent.
\end{thm}

More generally, the number of distinct orbits of $Aut(X,G,\Phi)$ in $X$, or, alternatively, the number of the classes of equivalent chains in $\fG(\Phi)$, measures the degree of non-homogeneity of the action $(X,G,\Phi)$. The number of  equivalence classes of group  chains in $\fG(\Phi)$ is invariant under topological conjugacy of dynamical systems. It is then a natural question to look for algebraic or geometric invariants which determine the number of the classes of equivalent chains in $\fG(\Phi)$.

The following  remark   allows us to restrict to the study of conjugate chains for a given chain $\{G_i\}_{i \geq 0}$ to those in a standard form.

\begin{remark}\label{representative}
{\rm
Let $(X,G,\Phi)$ be a minimal equicontinuous dynamical system, and let $\{G_i\}_{i \geq 0}$ be a group chain with basepoint $x$. Let $y \in X$, and consider the equivalence class of all group chains with basepoint $y$. Then this equivalence class has a representative of the form $\{g_iG_ig_i^{-1}\}_{i \geq 0}$, where $g_iG_i = g_jG_i$ for all $j \geq i$.
}
\end{remark}

 Finally, we will need  the following theorem, proved in Fokkink and Oversteegen \cite{FO2002}.

Recall that $\ds {\rm core}_G \, G_i = \bigcap_{h \in G} hG_ih^{-1}$ is a maximal subgroup of $G_i$ normal in $G$. If $G_i \supset G_{i+1}$, then ${\rm core}_G  \, G_i \supset {\rm core}_G  \, G_{i+1}$, and so $\{{\rm core}_G \, G_i\}_{i \geq 0}$ is a chain of subgroups of finite index in $G$.

\begin{thm}\cite{FO2002}\label{charact-groupchains}
Let $\{G_i\}_{i  \geq 0}$ be a group chain in $G$. Then the following is equivalent:
\begin{enumerate}
\item The chain $\{G_i\}_{i \geq 0}$ is equivalent to a regular  group chain $\{N_i\}_{i \geq 0}$, where $N_0 = G$ and $N_i$ is a normal subgroup of $G$ for all $i > 0$.
\item The chain $\{G_i\}_{i \geq 0}$ is equivalent to the chain $\{{\rm core}_G  \, G_i\}_{i \geq 0}$.
\item Every chain $\{g_i G_i g_i^{-1}\}_{i \geq 0}$, where $g_jG_i = g_iG_i$ for all $j \geq i$, is equivalent to $\{G_i\}_{i \geq 0}$.
\end{enumerate}
If any of these conditions are satisfied, then by    Theorem~\ref{thm-rephr}, it follows  that the associated Cantor minimal system $(G_{\infty}, G)$ is homogeneous.
\end{thm}
 
 The equivalences of Theorem~\ref{charact-groupchains} allow us to make precise the comparison of the notions of weakly normal in Definition~\ref{def-weaklynormal} and virtually regular in Definition~\ref{def-virtually-regular} for group chains. 
 
\begin{remark}\label{rmk-comparison}
{\rm
The group chain $\{G_i\}_{i \geq 0}$ is weakly normal if there exists an index $i_0 \geq 0$ such that for the truncated group chain
$\{G_i\}_{i \geq i_0}$,  each conjugate equivalent chain in $G_{i_0}$ is equivalent to $\{G_i\}_{i \geq i_0}$.  By part (1) of Theorem~\ref{charact-groupchains}, we then have that the group chain $\{G_i\}_{i \geq i_0}$ is equivalent to a group chain $\{N_i\}_{i \geq i_0}$, where $N_0 = G_{i_0}$ and $N_i$ is a normal subgroup of $G_{i_0}$ for all $i > i_0$.

The group chain $\{G_i\}_{i \geq 0}$ is virtually regular   if there is a   subgroup $G_0' \subset G_0$ of finite index   such that for the restricted group  chain  $\{G_i' = G_i \cap G_0'\}_{i \geq 0}$,
each conjugate equivalent chain in $G_0'$ is equivalent to $\{G_i'\}_{i \geq 0}$.  
By part (1) of Theorem~\ref{charact-groupchains}, we then have that the group chain $\{G_i'\}_{i \geq 0}$ is equivalent to a chain $\{N_i'\}_{i \geq 0}$, where $N_0' = G_0'$ and $N_i'$ is a normal subgroup of $G_0'$ for all $i > 0$.

Thus, the difference between the two notions is that,  for a weakly normal chain $\{G_i\}_{i \geq 0}$, each   group   $N_i$ for $i  \geq i_0$ is normal in $G_{i_0}$, while for a virtually regular group chain,  each   group   $N_i'$ for $i  \geq i$ is normal in $G_0'$, but the groups $N_i'$ need not be normal in $G_{i_0}$ for any $i_0 \geq 0$.

}
\end{remark}

\section{The discriminant group}\label{sec-quotient}

 In this section, we introduce the core and discriminant groups associated to a group chain. We then show that these  groups, defined using inverse limits, provide    representations for the Ellis group $\oG$  and the isotropy subgroup $\oG_x$ at $x \in X$, as introduced in Sections~\ref{sec-intro} and \ref{sec-ellis}. 

Let $(X,G,\Phi)$ be an equicontinuous Cantor minimal system, and let $\{G_i\}_{i \geq 0}$ be a group chain with partitions $\{\cP_i\}_{i \geq 0}$ and basepoint $x \in X$.  For each subgroup $G_i$, define its   core   
\begin{equation}\label{eq-core}
C_i = {\ds  {\rm core}_G   \, G_i = \bigcap_{g \in G} gG_ig^{-1}} ~ \subset ~ G_i
\end{equation}
  which is the maximal subgroup of $G_i$, normal in $G$. Then the quotient $G/C_i$ is a finite group, 
    and the collection  $\{C_{i} \}_{i \geq 1}$ forms a descending chain of normal subgroups of $G$. 
    
  \begin{lemma}\label{c-infty}
  There are well-defined homomorphisms of finite groups
      \begin{align}\label{eq-corebonding} \delta^{i+1}_i \colon   G/C_{i+1} \to G/C_{i} \ ,
      \end{align}
for which the resulting   inverse limit space   is a profinite group
  \begin{align} \label{cinfty-define} 
  C_{\infty} = \lim_{\longleftarrow} \, \left\{\delta^{i+1}_i \colon   G/C_{i+1} \to G/C_{i}   \right\} \ . 
   \end{align}

  \end{lemma}  
  
  \proof Observe that $C_{i} $ is the intersection of a finite number of subgroups of finite index, and thus $G/C_{i} $ is finite. The chain $\{C_i\}_{i \geq 0}$ is nested, and the maps \eqref{eq-corebonding} are induced by coset inclusion.
  
  Since $C_i$ are normal, the maps \eqref{eq-corebonding} are group homomorphisms. 
It follows from results in \cite{RZ2000} that the inverse limit space   $C_{\infty}$ defined in   \eqref{cinfty-define} is a profinite group. 
\endproof

Now consider the quotient $D_i = G_i/C_{i} $. Since $C_{i} $ is a normal subgroup in $G$, and hence also normal in $G_i$, the quotient $D_i$ is a group. This group is naturally identified with a subgroup   $\widetilde{D}_i \subset G/C_{i} $.

\begin{lemma}\label{h-kernels}
There are group homomorphisms $\theta^i_j \colon   D_i \to D_j$, such that $\theta^i_j \circ \theta^k_i = \theta^k_j$. Thus, we obtain  a   profinite group
\begin{equation}\label{eq-discriminant}
\cD_x = \lim_{\longleftarrow}\, \left \{\theta^i_j \colon  D_i \to D_j \right\} \ .
\end{equation}
\end{lemma}

\proof
Denote by $q_i \colon   G/C_{i} \to \left( G/C_{i} \right)/D_i \cong G/G_i$, then we have the following diagram:
\begin{align*} \xymatrix{ G/C_{1} \ar[d]^{q_1} & \ar[l] G/ C_{2} \ar[d]^{q_2}& \ar[l] G/ C_{3} \ar[d]^{q_3} & \ar[l] \cdots \\
G/G_1 & \ar[l] G/G_2  & \ar[l] G/G_3 & \ar[l]\cdots } \end{align*}

Since $D_i \subset G/C_{i} $, the restriction of $\delta^{i+1}_i \colon  G/C_{i+1}  \to G/C_{i} $ gives a group homomorphism 
  $$\theta^{i+1}_{i} \colon   D_{i+1} \to G/C_{i}  \ . $$ 
Since $G_{i+1} \subset G_i$, and $\theta^{i+1}_i$ is the inclusion of cosets, we have that $\theta^{i+1}_i(D_{i+1}) \subseteq G_i/C_{i} = D_i $. By a similar argument one checks that $\theta^{i+1}_{i-1} = \theta^i_{i-1} \circ \theta^{i+1}_i$, and then by \cite[Propositions 1.1.1 and 1.1.4]{RZ2000} the inverse limit $\cD_x$ exists, 
and is a  profinite group.
\endproof

\begin{defn}\label{def-discriminant}
Let $(X,G,\Phi)$ be an equicontinuous Cantor minimal system, and $\{G_i\}_{i \geq 0}$ be a group chain with partitions $\{\cP_i\}_{i \geq 0}$ and basepoint $x \in X$.  The profinite group $\cD_x$ in \eqref{eq-discriminant} is called the \emph{discriminant group} of $(X,G,\Phi)$ at $x \in X$.
\end{defn}




We next establish the  relation between the discriminant group and    the Ellis group for the   Cantor minimal system $(G_{\infty}, G)$ associated to the group chain $\{G_i\}_{i \geq 0}$.
The following result is Theorem~\ref{thm-quotientspace-1} of the Introduction.

\begin{thm}\label{thm-quotientspace}
Let    $(X,G,\Phi)$ be an   equicontinuous minimal action, and 
let $\{\cP_i\}_{i \geq 0}$ be a sequence of partitions of $X$ chosen as in Proposition~\ref{prop-AFpres} with basepoint $x \in X$,  and let $\{G_i\}_{i \geq 0}$ be the associated group chain.
Then there is a natural isomorphism of topological groups $\whTheta \colon  \oG \cong C_{\infty}$ such that the restriction $\whTheta \colon \oG_x \cong \cD_x$.
\end{thm}
\proof Recall that the factor maps  $\phi_i = \tilde{\phi}_i \circ \jmath_i  \colon   X \to G/G_i$ defined in \eqref{eq-finiterep} are $G$-equivariant. Then the left action induces    
  representations $\Theta_i \colon G \to {\rm Perm}(G/G_i)$. Let $x_i$ denote the coset of the identity in $G/G_i$, then the image $\Theta_i(G_i)$ consists of the permutations which fix $x_i$.  That is, $G_i$ is the isotropy group of $x_i$ for the action $\Theta_i$.
  
  The kernel ${\rm ker}(\Theta_i) \subset G$ of the map $\Theta_i$ consist of all elements in $G$ which fix every element of $G/G_i$. As the action of $G$ on $G/G_i$ is transitive,  ${\rm ker}(\Theta_i)$ is the intersection of the conjugates of $G_i$ and thus equals the group $C_i$.  Thus, the representation $\Theta_i$ factors through the group $G/C_i$ which we identify with its image to give
\begin{equation}\label{eq-coordinates}
\whTheta_i \colon G \to G/C_i \subset {\rm Perm}(G/G_i) \ .
\end{equation}
Since $G/C_i$ is a finite set, hence has the discrete topology, the map $\whTheta_i$ extends to a continuous homomorphism
$\ds \whTheta_i \colon \whE(X,G,\Phi) \to  G/C_i$, where $\whE(X,G,\Phi) \subset X^X = {\it Maps}(X,X)$ is defined in Section~\ref{sec-ellis}. Recall from \eqref{cinfty-define} that $\ds C_{\infty} = \lim_{\longleftarrow}\, \left\{\delta^{i+1}_i \colon   G/C_{i+1} \to G/C_{i} \right\}$.
It follows that the product of the collection of maps $\{\whTheta_i \mid i \geq 1\}$
 defines a continuous homomorphism $\whTheta \colon \whE(X,G,\Phi) \to C_{\infty}$.
We claim that $\whTheta$ is an isomorphism of topological groups. 

The map $\whTheta$ is surjective as each map $\whTheta_i$ is surjective. 

Suppose that $\whphi \in \whE(X,G,\Phi) \subset X^X$ is not the identity. Then for some $x \in X$ we have $ \whphi(x) = y \ne x$. Then $\e = d(x,y) > 0$, and choose $i$ sufficiently large so that $2^i < \epsilon/2$. Then by condition (2) of Proposition~\ref{prop-AFpres} the points $x,y$ lie in disjoint elements of the partition $\cP_i$. It follows that 
$\whTheta_i(\whphi)$ is not the identity. Thus, $\whTheta(\whphi)$ is not the identity, and $\whTheta$ is injective. Since $\whE(X,G,\Phi)$ is identified with $\oG$, we have an isomorphism between $\oG$ and $C_{\infty}$.

 Finally, we show that  $\whE(X,G,\Phi)_x= \oG_x$ is identified with $\cD_x$.
 The group $G_i$ is the isotropy of $\Phi$ at $x_i$, so $G_i \subset \whE(X,G,\Phi)_x$ and thus has image   $\whTheta_i(G_i) = G_i/C_i \equiv D_i \subset G/C_i$ where the notation $D_i$ was introduced in Lemma~\ref{h-kernels}. Again, as $D_i$ is discrete, we obtain a continuous homomorphism  
 $\whTheta_i \colon \whE(X,G,\Phi)_x \to D_i$ for each $i \geq 1$. From the definition   of $\cD_x$ in \eqref{eq-discriminant}, we thus obtain a continuous homomorphism   $\whTheta \colon  \whE(X,G,\Phi)_x \to \cD_x$.
 Then as before, this map is a topological isomorphism.
  \endproof

The \emph{tree diagram}   associated to a collection of partitions  $\{\cP_i\}_{i \geq 0}$ of $X$ with basepoint $x$ is a tree $\cT$,  whose vertices      at level $i$ correspond to the clopen sets in $\cP_i$. For each $i \geq 0$ the set $V_i \in \cP_i$ is determined by the condition $x \in V_i$, so the vertices of $\cT$ are identified with the points in the coset $G/G_i$ via its action on the clopen sets. The edges of $\cT$ are defined by the inclusions of clopen sets of $\cP_i$ in $\cP_{i-1}$, which corresponds to the quotient maps $G/G_i \to G/G_{i-1}$.   This construction is discussed further in \cite[Chapter 2]{Dyer2015}. Each point $y \in X$ is assigned a unique path in the tree $\cT$ associated to the nested chain of partition clopen sets containing $y$. In this sense, the tree   $\cT$ associated to $\{\cP_i\}_{i \geq 0}$   can be considered as a ``coordinate system'' for the Cantor set $X$. This construction is related to the Bratteli diagram for $\mZ$-actions introduced   by  Ver{\v{s}}ik  in \cite{Vershik1982}, and   discussed in  \cite{HPS1992} for example.

The map $\whTheta \colon \oG \to C_{\infty}$ constructed in Theorem~\ref{thm-quotientspace} represents the action of $\oG$ on $X$ as automorphisms of the tree $\cT$.   The mapping $\whTheta_i \colon G \to G/C_i$ in \eqref{eq-coordinates} gives a finite approximation to the action of $G$ on $X$, as permutations of the finite sets $G/G_i$ which are the vertices of the tree $\cT$ at level $i$.  Thus, Theorem~\ref{thm-quotientspace} can be interpreted as providing  a representation  for the   Ellis group $\oG$ of $(X,G,\Phi)$ in   ``tree coordinates'' on $X$. In this representation, the isotropy group $\oG_x$ for $x \in X$ is mapped to a subgroup of the automorphisms of  $\cT$ which fix the path corresponding to the basepoint $x$. The automorphisms induced by the action of $\oG_x$ on the vertices $G/G_i$ of $\cT$ at level $i$ correspond to the action of the quotient group $G_i/C_i$.

\section{Properties of the discriminant group}\label{sec-properties}

In this section, we show three basic properties of the discriminant group, which are used to investigate   equicontinuous   Cantor minimal systems. First,  we consider how the discriminant group $\cD_x$ depends on the choices which are made in its definition. Second, we show that the discriminant group is ``completely not normal''. Finally, we introduce the notion of ``normal form'' for group chains, with respect to the discriminant group.
 
 Let $(X,G,\Phi)$ be an equicontinuous Cantor minimal system.
  Let    $\{G_i\}_{i \geq 0}$  be a group chain associated to a  collection of partitions $\{\cP_i\}_{i \geq 0}$ with basepoint    $x \in X$, and let    $\cD_x$ be its discriminant group, as defined by \eqref{eq-discriminant}.
Consider  another group chain  $\{H_i\}_{i \geq 0}$ for $(X,G,\Phi)$ associated  with partitions $\{\cQ_i\}_{i \geq 0}$ with basepoint  $y \in X$, and   let    $\cD_y'$ be its discriminant group, as defined by \eqref{eq-discriminant}.

\begin{prop}\label{prop-isomorphism}
There is an isomorphism $\cD_x \cong \cD_y'$. Thus, the cardinality of $\cD_x$ is an invariant of the homeomorphism class of the equicontinuous Cantor minimal system $(X,G,\Phi)$.
\end{prop}
\proof
Let $\oG$ be the Ellis group of the action $(X,G,\Phi)$. Then by Theorem~\ref{thm-quotientspace}, the core group  $C_{\infty}$ of $\{G_i\}_{i \geq 0}$  and the core group $C_{\infty}'$ of  $\{H_i\}_{i \geq 0}$ are both topologically   isomorphic to $\oG$, and thus $C_{\infty} \cong  C_{\infty}'$. By the identity \eqref{eq-profiniteconj}, the isotropy groups $\oG_x$ and $\oG_y$ are conjugate in $\oG$. Then   by Theorem~\ref{thm-quotientspace} we have that $\cD_x \cong \oG_x \cong \oG_y \cong \cD_x'$.
\endproof

The isomorphism of the groups $\cD_x$ and $\cD_y'$ in Proposition~\ref{prop-isomorphism} is abstract, as it invokes their relation with subgroups of the Ellis group $\oG$.  
The  thesis of the first author \cite[Theorem~7.5]{Dyer2015} gives a more precise form of  this isomorphism in the case where the point $y \in G \cdot x$,  using the calculus of group chains. Given the Cantor minimal system $(X,G,\Phi)$,  let $\{G_i\}_{i \geq 0}$ and $\{H_i\}_{i \geq 0}$ be group chains associated with the action, with common basepoint $x$. Then by 
 Theorem~\ref{thm-equiv-rel}, the group chains are equivalent, and the equivalence induces a canonical isomorphism of their core chains $C_{\infty} \cong  C_{\infty}'$. Simple ``diagram chasing'' shows that this isomorphism restricts to an isomorphism $\cD_x$ to $\cD_x'$. This direct approach extends to the case where $\{H_i\}_{i \geq 0}$ is a group chain with basepoint $y \in   G \cdot x$.

 In the case where  the basepoint $y$ of the chain $\{H_i\}_{i \geq 0}$  is in the orbit of $x \in X$ under the action of $Aut(X,G,\Phi)$, then we have the following more precise  statement, which is related to the considerations in Section~\ref{sec-kernel}.

\begin{cor}\label{cor-constant-auto}
Let $(X,G,\Phi)$ be an equicontinuous Cantor minimal system  
 with associated group chain $\{G_i\}_{i \geq 0}$ at $x$, and let $\{g_iG_ig_i^{-1}\}_{i \geq 0}$ be a group chain at $y$. Let $\cD_x$ and $\cD_y$ be the discriminant groups at $x$ and $y$ respectively. If there exists $h \in Aut(X,G,\Phi)$ such that $h(x) = y$, then $\cD_x = \cD_y$. That is, the groups are identical as subgroups of $G_{\infty}$.
\end{cor}

\proof 
As shown in the proof of Proposition~\ref{prop-auto=normalizer}, the assumption that $x,y$ are in the same orbit of ${Aut}(X,G,\Phi)$  implies that there exists $(g_i)    \in \whcN(\oG_x)$, where $\whcN(\oG_x)$ is the normalizer of $\oG_x$ in $\oG$, so that 
$(g_i) \cdot x = y$. Then by  \eqref{eq-profiniteconj}
$$\oG_y =   \oG_y = (g_i) \cdot \oG_x \cdot (g_i)^{-1} = \oG_x.$$
Since the limit core $C_{\infty}$ is independent of the choice of a point, this implies $\cD_x = \cD_y$.
\endproof

Next, let  $\{G_i\}_{i \geq 0}$ be a group chain with discriminant group $\cD_x$.   The countable subgroup $G \subset C_{\infty}$ can be thought of as the ``rational points'' in $C_{\infty}$, and we introduce also the ``rational'' core group
\begin{equation}
\cC_x = \bigcap_{k \in G} ~ k \  \cD_x \ k^{-1} \subset \cD_x \ .
\end{equation}

\begin{prop}\label{prop-discrisnotnormal-1}
  The rational core subgroup $\cC_x$ is trivial.
\end{prop}

\proof We have the definitions 
\begin{eqnarray*}
C_{\infty} & = & \lim_{\longleftarrow} \{\delta^i_{i-1}:G/C_i \to G/C_{i-1}\}\\
\cD_x & = & \lim_{\longleftarrow} \{\delta^i_{i-1}:G_i/C_i \to G_{i-1}/C_{i-1}\}
\end{eqnarray*}
  where  $C_i$ is the maximal normal subgroup of $G_i$.    For $i \geq 1$,
 denote by $\delta_i \colon C_{\infty} \to G/C_i$ the projection map, which commutes with both the left and the right actions of $G$.  
Since each $\delta_i$ is a  group homomorphism, the image $\delta_i(\cD_x)$ is a subgroup of $G_i/C_i \subset G/C_i$ and we have
$$
\delta_i(\cC_x) = \delta_i(\bigcap_{k \in G} ~ k \  \cD_x \ k^{-1}) =  \bigcap_{k \in G} ~ k \  \delta_i(\cD_x) \ k^{-1}
\subset   \bigcap_{k \in G} ~ k \  (G_i/C_i) \ k^{-1} \subset C_i/C_i  
$$
which is the trivial group. Thus, $\delta_i(\cC_x)$ is the trivial group for all $i \geq 1$, and hence $\cC_x$ must be the trivial group.
\endproof
 
\begin{cor} \label{cor-maximalsubgroup}
The maximal normal subgroup of $\cD_x$ in $C_{\infty}$ is trivial. 
\end{cor}
\proof 
Observe that
 $$  {\rm core}_{C_{\infty}}  \cD_x =  \bigcap_{(g_i) \in C_{\infty}} ~ (g_i) \  \cD_x \ (g_i)^{-1} ~ \subset ~  \bigcap_{k \in G} ~ k \  \cD_x \ k^{-1}  = \cC_x$$
which implies the claim. 
\endproof
 
 \bigskip
 
 We   obtain the following result, which implies Theorem \ref{cor-discrisnotnormal} and Corollary \ref{cor-3-1} of the Introduction.
 
 \begin{thm}\label{thm-normaldiscr}
Let $(X,G,\Phi)$ be an equicontinuous Cantor minimal system, and let $\oG$ be its Ellis group. Then for $x \in X$, the isotropy group $\oG_x$ is a normal subgroup of $\oG$ if and only if $\oG_x$ is trivial. 
\end{thm}

\proof 
If $\oG_x$ is trivial, the claim is immediate, as its normalizer is $\oG$. 

Conversely, suppose $\oG_x$ is a normal subgroup of $\oG$. 
 Let    $\{\cP_i\}_{i \geq 0}$ be a sequence of partitions   with basepoint $x \in X$, and let  $\{G_i\}_{i \geq 0}$ be the associated  group chain. Then $\oG_x$  is isomorphic to the discriminant group $\cD_x$ of $\{G_i\}_{i \geq 0}$  by Theorem~\ref{thm-quotientspace}. By Corollary~\ref{cor-maximalsubgroup}, the maximal normal subgroup of $\cD_x$ is trivial, hence the same hold for $\oG_x$. 
\endproof

We conclude this section with a  useful technical result, that for  any Cantor minimal system  $(X,G,\Phi)$, there exists a choice of a group chain $\{G_i\}_{i \geq 0}$ associated to the action which is in \emph{normal form}:

\begin{defn}\label{def-normalform}
 Let $(X,G,\Phi)$ be a minimal equicontinuous action, and let $\{G_i\}_{i \geq 0}$ be a group chain with partitions $\{\cP_i\}_{i \geq 0}$ and basepoint $x \in X$. We say that  $\{G_i\}_{i \geq 0}$ is \emph{in the normal form} if   the projections $\cD_x \to G_i/C_i$ are onto for every $i \geq 0$.
\end{defn}

 The next result shows that a representing group chain in normal form always exists.

 \begin{prop}\label{prop-specialchoice}
 Let $(X,G,\Phi)$ be a minimal equicontinuous action. Then there exists a collection of partitions $\{\cP_i\}$ of $X$, such that the associated group chain $\{G_i\}_{i \geq 0}$ with basepoint $x \in X$ is in normal form.
  \end{prop}
 
 \proof Let $\{H_i\}_{i \geq 0}$ be a group chain with partitions $\{\cQ_i\}_{i \geq 0}$ and basepoint $x \in X$, and let 
  ${\ds C_i = \bigcap_{g \in G}gH_ig^{-1}}$ be the core of $H_i$. Let $\ds C_{\infty} = \lim_{\longleftarrow}\, \left\{G/C_i \to G/C_{i-1}\right\}$, and let  
  $$\cD_x' = \lim_{\longleftarrow}\, \left\{H_i/C_i \to H_{i-1}/C_{i-1} \right\} \subset C_{\infty}$$ 
  be the discriminant group at $x$. Let $D_i = \theta_i(\cD_x')$, where $\theta_i \colon   C_{\infty} \to H_i/C_i$ are the projections. 
  
  Since $C_i$ is a normal subgroup, the union $G_i = \bigcup \{a C_i \mid aC_i \in D_i\}$ is a subgroup of $G$. Since $\{C_i\}_{i \geq 0}$ is a nested chain, and the bonding maps $D_i \to D_{i-1}$ are coset inclusions, the collection $\{G_i\}_{i \geq 0}$ forms a nested chain of subgroups of $G$.
 
 We claim that the group chain $\{G_i\}_{i \geq 0}$ is equivalent to $\{H_i\}_{i \geq 0}$. First we note that by construction,  $G_i \subset H_i$ for every $i \geq 0$. 
 Next, we claim that for every $i \geq 0$ there exists $j \geq i$ such that $H_j \subset G_i$. Let $bC_i \in H_i/C_i$. If $b \notin G_i/C_i = D_i$, then there exists $j_b>i$ such that $ bC_i \notin \theta^{j_b}_i(H_{j_b}/C_{j_b})$. Since $H_i/C_i$ is finite, we can define
   \begin{align*} 
   j = \max\{j_b \mid bC_i \notin D_i\} \ . 
   \end{align*}
Then $\theta^j_i(H_j/C_j) = D_i$. By possibly restricting to subsequences, we can assume that 
   $$H_1 \supset G_1 \supset H_2 \supset G_2 \supset \cdots \ ,$$
 and so $\{G_i\}_{i \geq 0}$ and $\{H_i\}_{i \geq 0}$ are equivalent group chains.
 
 Then by Corollary \ref{cor-conjequiv} there exists a collection of partitions $\{\cP_i\}_{i \geq 0}$ such that $\{G_i\}_{i \geq 0}$ is the associated group chain with basepoint $x$. By construction, $C_i \subset G_i$, and since $C_i$ is maximal in $H_i$ and $G_i \subset H_i$, then $C_i$ is also maximal in $G_i$. So $C_i = {\rm core}_G  \, G_i$, and then $C_{\infty}$ is the limit core associated to $\{G_i\}_{i \geq 0}$. Let $\ds \cD_x = \lim_{\longleftarrow} \, \left\{G_i/C_i \to G_{i-1}/C_{i-1} \right\}$. Then by construction $\cD_x = \cD_x'$, and the projections $\cD_x \to G_i/C_i$ are onto.
 \endproof

  \section{Finite discriminant group}\label{sec-weakly-reg-dx}

In this section, we give a result which characterizes the  equicontinuous Cantor minimal systems with finite discriminant group $\cD_x$.

 \begin{thm}\label{thm-core-finite}
 Let $(X,G,\Phi)$ be an equicontinuous Cantor minimal system, and let $\{G_i\}_{i \geq 0}$ be a group chain with partitions $\{\cP_i\}_{i \geq 0}$ and basepoint $x \in X$.
 Then the discriminant group $\cD_x$ is finite if and only if $\{G_i\}_{i \geq 0}$ can be chosen so that  for all $i \geq 1$, we have
    \begin{align}\label{eq-kernel}G_{i} \cap C_{1} = C_{i}. \end{align}
 \end{thm}
 
 \proof 
Denote by $\delta_i \colon   \cD_x \to G_i/C_i$ the projections. 
Note that 
  \begin{align}\label{eq-kernelcond}(\delta^{i}_1)^{-1} (eC_1) = (G_i \cap C_1)/C_i \ . 
  \end{align}
 
 By Proposition~\ref{prop-specialchoice} we can assume that the projections $\delta_i \colon   \cD_x \to G_i/C_i$ are surjective for all $i \geq 0$.
  
If $\cD_x$ is finite, then there exists $m>0$ such that $card(\cD_x) = card(G_i/C_i) = card(D_i)$ for all $i \geq m$. By discarding a finite number of groups $G_1,\ldots, G_{m-1}$ from the sequence, we can assume that $m=1$. 
 
 Then the preimage $\delta_{1}^{-1}(eC_1)$ consists of exactly one element, and so every bonding map 
   $$\delta^i_1 \colon  G_i/C_i \to G_1/C_1$$ 
  is injective, that is, $(\delta^i_1)^{-1}(eC_1) = eC_i$. Then \eqref{eq-kernelcond} implies that $G_i \cap C_1 = C_i$ for all $i \geq 1$. 
 
 For the converse, suppose a group chain $\{G_i\}_{i \geq 0}$ with property \eqref{eq-kernel} exists. Since $G_i \cap C_1 = C_i$, the preimage $(\delta^i_1)^{-1}(eC_1) = e C_i$ for all $i \geq 1$, and so the preimage $\delta_i^{-1}(eC_1) \subset \cD_x$ is a single element. Since $\cD_x$ is a group, and $G_1/C_1$ is a finite group, this implies that $\cD_x$ is finite. 
 \endproof
 
 The proof of Theorem~\ref{thm-virtually-regular} of the introduction follows immediately by taking $N = C_1$.

 \begin{ex}\label{intersection-regular}
 {\rm
 Suppose the action $(X,G,\Phi)$ is homogeneous, that is, $Aut(X,G,\Phi)$ acts transitively on $X$. Then by Theorem \ref{charact-groupchains} the group chain $\{G_i\}_{i \geq 0}$ at $x$ can be chosen so that every $G_i$, $i \geq 0$, is a normal subgroup of $G$. In this case $G_i = C_i$, $i \geq 0$, and then $C_1 \cap G_i = C_i = G_i$.
 }
 \end{ex}

We next deduce the application  of Theorem~\ref{thm-core-finite}  to weak solenoids,   given in Section~\ref{sec-intro} as Corollary~\ref{cor-almostcover}, which we restate for convenience: 
 \begin{cor}\label{cor-normalcover}
 Let $M_0$ be a closed manifold, $x_0 \in M_0$ a fixed basepoint, and let $G = \pi_1(M_0 , x_0)$ be the fundamental group.
Let $f \colon \cS \to M_0$ be a weak solenoid, with  $X= f^{-1}(x_0)$ the Cantor fiber, and $x \in X$ a given basepoint. 
 Suppose the   monodromy action $(X,G,\Phi)$ on the fiber has   finite discriminant group $\cD_x$. Then there exists a finite-to-one covering map $p_0 \colon M_0' \to M_0$ such that the pullback solenoid $f' \colon  \cS' \to M_0'$ is homogeneous.
 \end{cor}
 
 \proof 
 
 If $(X,G,\Phi)$ is an action with finite discriminant, then by Theorem~\ref{thm-core-finite} the group chain $\{G_i\}_{i \geq 0}$ can be chosen so that $G_i \cap C_1 = C_i$, where $C_i$, $i \geq 1$, is the normal core of $G_i$. That is, for each $i \geq 0$ there is a closed manifold $M_i = \widetilde{M}_0/G_i$, where $\widetilde{M}_0$ is the universal cover of $M_0$, and by the standard argument we have maps $f^i_{i-1}: M_i \to M_{i-1}$, and $\ds \cS \cong \lim_{\longleftarrow} \, \left\{f^i_{i-1}:M_i \to M_{i-1} \right\}$.
 
Since $C_1$ is a normal subgroup of $G$ of finite index, there exists a finite-to-one covering $p_0: M_0' \to M_0$, with $\pi_1(M_0',x_0') = C_1$. Then each $f^i_{i-1}: M_i \to M_0$ pulls back to a map $g^i_{i-1}: M_i' \to M_0'$, where $\pi_1(M_i',x_i') = G_i \cap C_1 = C_i$ is a normal subgroup of $C_1$. 
Thus the weak solenoid $\ds \cS' = \lim_{\longleftarrow} \, \left\{M_i' \to M_{i-1}' \right\}$ is homogeneous.
 \endproof

Suppose that $\cS$ is a weak solenoid which defines the group chain  $\{G_i\}_{i \geq 0}$ as in Section~\ref{sec-intro}. Recall that $M_0$ denotes the base manifold for $\cS$. Then in the weakly normal case,  the covering of $M_0$ defined by the subgroup $G_{i_0} \subset G_0 = \pi_1(M_0, x_0)$ is simply $M_{i_0}$ which is one of the coverings used to define the inverse limit presentation of $\cS$, and the induced solenoid on $M_{i_0}$ is homogeneous. In contrast, if $\{G_i\}_{i \geq 0}$  is virtually regular, then for the covering $M_N$ of $M_0$ corresponding to the subgroup $N \subset G_0$, then   the induced solenoid on $M_0'$ as defined in Corollary~\ref{cor-normalcover} is homogeneous, but this need not imply that $\cS$ is homogeneous.

\section{Kernel of a group chain}\label{sec-kernel}

In this section we investigate another property of group chains, the kernel   of a group chain 
  $\{G_i\}_{i \geq 0}$   of finite index subgroups in a finitely generated group $G$. 

\begin{defn}\label{FO2002}
Let $(X,G,\Phi)$ be a minimal equicontinuous action, and let $\{G_i\}_{i \geq 0}$ be a group chain at $x \in X$.
The \emph{kernel} of the group chain $\{G_i\}_{i \geq 0}$ at $x$ is the intersection
  $\ds K_x = \bigcap_i G_i.$
\end{defn}

The following possibilities can be realized for the kernel of a group chain $\{G_i\}_{i \geq 0}$:
\begin{enumerate}
\item The kernel $K_x$ is trivial, that is, $K_x = \{e\}$, where $e$ is the identity element in $G$.  
\item The kernel $K_x$ is a non-trivial normal subgroup of $G$.
\item The kernel $K_x$ is a non-trivial non-normal subgroup of $G$.
\end{enumerate}

We note that if $\{H_i\}_{i \geq 0}$ is a group chain equivalent to $\{G_i\}_{i \geq 0}$, then   for every $i \geq 0$ there exists $j \geq i$ such that $G_j \subset H_i$ and $H_j \subset G_i$, and it follows that   
  $$K_x = \bigcap_i G_i = \bigcap_i H_i \ .$$
Thus,  the kernel of a group chain is invariant under the equivalence of group chains.  Example~\ref{ex-RT} below shows that the kernel of a group chain is not invariant under the conjugate equivalence of group chains. Consider also the maximal normal subgroup of the kernel, 
\begin{equation}\label{eq-kernelcore}
L_x = {\rm core}_G  \, K_x \ .
\end{equation}
Then $L_x$ is an invariant of the conjugate equivalence class, so in particular it is independent of $x$.  
It  can be identified with the kernel of the homomorphism $\Phi \colon G \to \oG$, as seen in the following:

\begin{lemma}\label{lem-kernel}
Let $(X,G,\Phi)$ be a minimal equicontinuous group action with group chain $\{G_i\}_{i \geq 0}$ at $x \in X$, let $K_x = \bigcap_i G_i$, and let $L_x = {\rm core}_G  \, K_x$. 
Then the action $G/L_x$ on $X$ is effective. If $\cD_x$ is trivial, then the action $G/L_x$ on $X$ is free, and so $K_x = L_x$.
\end{lemma}

\proof To show that the action of $G/L_x$ on $X$ is effective, we need to show that if $h \notin L_x$, then there exists $y \in X$ such that $h \cdot y \ne y$. Let $y \in X$, then by Remark \ref{representative} there is a group chain $\{g_iG_ig_i^{-1}\}_{i \geq 0}$ at $y$ with $K_y = \bigcap_i g_iG_ig_i^{-1}$, and if $h \cdot y = y$, then $h \in K_y$. Then if $h$ fixes every point in $X$, this implies that $h \in \bigcap_{g \in G} gG_ig^{-1}$ for all $G_i$ and so $h \in L_x$. Therefore, if $h \notin L_x$, there must be $y \in X$ such that $h \cdot y \ne y$.

Now let $\cD_x$ be trivial, and let $g \in K_x$, that is, $g \cdot x = x$. For every $y \in X$ there is $h \in Aut(X,G,\Phi)$ such that $h(x) = y$, and so $g \cdot y = g \cdot h(x) = h(g \cdot x) = h(x) =y$. Then $g \in L_x$, and the action of $G/L_x$ on $X$ is free.
\endproof

Lemma~\ref{lem-kernel} states that the kernel of a group chain defining a homogeneous action is always a normal subgroup of $G$. A more interesting problem is to investigate  the relationship between the kernel of a group chain and the discriminant group in the case when the discriminant group is non-trivial.
It turns out that sometimes one can represent elements of the discriminant group using the elements in the kernel of a group chain.

\begin{prop}\label{prop-specialkernel}
Let $(X,G,\Phi)$ be an equicontinuous Cantor minimal system, and let $\{G_i\}_{i \geq 0}$ be a group chain at $x \in X$. Suppose $\{G_i\}_{i \geq 0}$ is in the normal form, that is, for any $i \geq 1$ the bonding maps $G_i/C_i \to G_{i-1}/C_{i-1}$ are surjective. Suppose that for every $i \geq 0$ we have
  \begin{align}\label{cond-kernel}G_i = K_xC_i,\end{align}
where ${\ds C_i = \bigcap_{g \in G} gG_ig^{-1}}$. If $K_x$ is a finite subgroup of $G$, then the discriminant group $\cD_x$ is finite.
\end{prop}

 \proof Since the group chain $\{G_i\}_{i \geq 0}$ is in the normal form, the bonding maps $G_i/C_i \to G_{i-1}/C_{i-1}$ are surjective, and so, since $G_i = K_x C_i$, we have the following commutative diagram.
   \begin{align}\label{eq-diagramkernel} \xymatrix{G_{i-1}/C_{i-1}  \ar[r]^{\cong} & K_xC_{i-1}/C_{i-1} \ar[r]^{\cong} & K_x/K_x \cap C_{i-1} \\
   G_{i}/C_{i} \ar[u]^{\delta^i_{i-1}} \ar[r]^{\cong} & K_xC_{i}/C_{i} \ar[u]^{\delta^i_{i-1}} \ar[r]^{\cong} & K_x/K_x \cap C_{i} \ar@{-->}[u]_{\kappa^i_{i-1}}
   }  \end{align}
This diagram defines the maps $\kappa^i_{i-1}: K_x / K_x \cap C_i \to K_x / K_x \cap C_{i-1}$, which, by commutativity of the diagram, must be surjective group homomorphisms, and so there is the profinite group $\ds \mathcal{K}_x = \lim_{\longleftarrow} \, \left\{K_x / K_x \cap C_i \to K_x / K_x \cap C_{i-1} \right\}$, and the group isomorphism $\cD_x \to \mathcal{K}_x$.

If $K_x$ is a finite group, then there exists $k \geq 0$ such that for all $i >k$ the maps $\kappa^i_{i-1}$ are the identity maps, and the group $\mathcal{K}_x$ is finite. Then the discriminant $\cD_x$ is a finite group.
 \endproof

   Proposition \ref{prop-specialkernel} is by no means exhaustive of all the possibilities when the discriminant or the kernel are non-trivial, and describes only one special situation.  We note that the property \eqref{cond-kernel} need not be invariant under the change of the basepoint $x$, and there are examples of both alternatives,  as shown in Section \ref{sec-examples}. 

 Recall that the definition of a stable kernel was given in Definition~\ref{def-stable}. We then note the above discussion implies the following basic result.  
\begin{prop}\label{prop-unstable}
Let $\{G_i\}_{i \geq 0}$ be a group chain, and suppose the corresponding discriminant group $\cD_x$ is nontrivial. If the kernel $K_x$ is a normal subgroup of $G$, then $\cD_x$ is unstable. 
\end{prop}
 \proof
The assumption implies that  $K_x \subset C_i$ for all $i \geq 0$, and so no element of $\cD_x$ can be represented by an element of $K_x$.
\endproof

 Finally, we give the promised examples of group chains with both stable and unstable kernels. 
The first example is similar to the example by Rogers and Tollefson in \cite{RT1972}, as described by Fokkink and Oversteegen in \cite[page 3750]{FO2002}; we thank the referee for bringing this example to our attention.

\begin{ex}\label{ex-RT}
{\rm
Let $G = \{ (a,b) \mid b^2 = e, bab = a^{-1} \}$ be infinite dihedral group, and consider the subgroups $G_i = \{(a^{2^i},b) \mid a,b \in \mZ\}$, for $i \geq 0$. Notice that the kernel $K_x = \bigcap G_i = \{e, b\}$, which is a finite subgroup isomorphic to $\mZ_2$. Here $x = (eG_i)$.

A simple computation shows that $|G:G_i| = 2^i$, and every coset of $G/G_i$ is represented by $a^k$, $0 \leq k < 2^{i}$, that is, there is a bijection $G/G_i \to \mZ/2^i \mZ$. However, the quotient $G/G_i$ is not a group, that is, cosets of $G/G_i$ do not act on each other. Indeed, $b \in eG_i$, and $baG_i = a^{-1}bG_i = a^{-1}G_i$, while the action of the identity fixes each coset.

The maximal normal subgroup of $G_i$ in $G$ is $C_i = \{a^{2^i}\}$, and so $G_i/C_i = \{eC_i, bC_i\}$ is a finite subgroup. Note that $G_i = K_x C_i$, that is, the hypothesis of Proposition \ref{prop-specialkernel} is satisfied. It follows that $\cD_x$ is finite and isomorphic to $K_x \cong \mZ_2$. The non-trivial element of $\cD_x$ is $(bC_i)$, with $b \in K_x$.

 In \cite{FO2002}, Fokkink and Oversteegen computed the kernels of chains associated to points in $G_{\infty}$ for the Rogers and Tollefson example. A very similar argument can be made to compute the kernels of group chains here, showing that if ${\ds y \in \lim_{\longleftarrow}\{G/G_i \to G/G_{i-1}\}}$, then either $K_y = \langle a^{2n} b \rangle $ for some integer $n\geq 0$, or $K_y = \{e\}$. 
 
 This shows that the kernel of a group chain is not invariant under the conjugate equivalence of group chains, since, depending on the choice of a basepoint, it may be trivial or non-trivial. Also, this shows that condition \eqref{cond-kernel} in Proposition \ref{prop-specialkernel} is not invariant under the change of a basepoint $x$, and so there are choices of a group chain for this example such that  an element of the discriminant cannot be represented by an element of the kernel of a group chain.
 
 Note that for the subgroup $G_0' = \{ (a,0)  \mid a \in \mZ \}$ of index $2$ in $G$, the chain $\{G_i' = G_i \cap G_0'\}$ is normal. Hence the group chain $\{G_i\}_{i \geq 0}$ is virtually regular, but is not weakly normal. 
  }
 \end{ex}

Next, we consider  examples where the elements of the discriminant group can be represented by elements of the kernel for any basepoint in $G_{\infty}$.

\begin{ex}\label{finitediscr}
{\rm
Let $\G$ be a finitely generated group, and $\{\G_i\}_{i \geq 0}$ be a normal group chain  in $\G$ with kernel  
   $\ds \G_x = \bigcap_i \G_i$. 
 Let $H$ be a finite simple group, and let $K \subset H$ be a non-trivial subgroup. Since $H$ is simple, $K$ is not normal in $H$.

Let $G = H \times \G$, and $G_i = K \times \G_i$, $i \geq 0$. We show that the discriminant group $\cD_x$ of this action, where $x = (eG_i)$, is finite.

First note that $G/G_i = (H \times \G)/(K \times \G_i) = H/K \times \G/\G_i$. Denote $\ds \G_{\infty} = \lim_{\longleftarrow} \,\left\{\G/\G_i \to \G/\G_{i-1} \right\}$. Then $G_{\infty} = H/K \times \G_{\infty}$. To compute the discriminant group, we have to compute the normal cores $C_i$ of $G_i$. Since $G$ is a product,  and $H$ is simple, then we have
  $$C_i = \bigcap_{g \in G} g G_ig^{-1} = \{e\} \times \G_i.$$
Then $G_i/C_i = K \times \G_i/\G_i$, and   $\delta^i_{i-1} \colon  G_i/C_i \to G_{i-1}/C_{i-1}$    is an injective map of spaces, isomorphic to $K$. Therefore, $\cD_x = \{(b,eC_i) \mid b \in K\}\cong K$ is finite. 

Now let $\{H_i =g_iG_ig_i^{-1}\}_{i \geq 0}$ be a chain of conjugate subgroups at $y$, where $g_i=(c_i,\gamma_i)$. Note that the projections $G/G_i \to G/G_{i-1}$ are the identities when restricted to $H/K$, so one can write $g_i = (c, \gamma_i)$, and then, since $\G_i$ are normal subgroups, $g_i G_i g_i^{-1} =  cKc^{-1} \times \G_i$.   Then the kernel of the chain $\{H_i\}$ is given by $K_y = \bigcap H_i = cKc^{-1} \times \G_x$. It is easy to see that the discriminant $\cD_y \cong cKc^{-1} \cong K$, and, in particular, every element in $\cD_y$ can be realized using an element in $K_y$. More precisely, $\cD_y = \{(b,eC_i) \mid b \in cKc^{-1}\}$.

}
\end{ex}

\begin{remark}
 {\rm Thus, the discriminant group in Example~\ref{ex-RT}   is unstable, and the discriminant group in Example \ref{finitediscr} is stable.
}
 \end{remark}

\section{Examples}\label{sec-examples}

In this final section, we give more examples of actions of non-abelian groups, illustrating the ideas discussed   in this paper. Our examples include the actions of the Heisenberg group, and of the generalized dihedral group.

\subsection{Heisenberg group}
Let $\cH$ be the discrete Heisenberg group, presented in the form $\cH = (\mZ^3, *)$ with the group operation $*$ given by $(x,y,z)*(x',y',z')=(x+x',y+y',z+z'+xy')$. This operation is standard addition in the first two coordinates, and addition with a twist in the last coordinate. We think about $H$ as $\mZ^2 \times \mZ$, where the $\mZ^2$ part is abelian, and the $\mZ$ part is not. 

As in \cite{LSU2014}, consider the subgroups of $\cH$ which can be written in the form $\Gamma=M\mZ^2\times m\mZ$ where $M = \left( \begin{array}{cc} i & j \\ k & l\end{array}\right)$ is a 2-by-2 matrix with non-negative integer entries and $m>0$ is an integer. Then $\gamma \in \Gamma$ is of the form $\gamma = (ix+jy, kx+ly, mz)$ for some $x,y,z \in \mZ$. A straightforward computation gives the following lemma.

\begin{lemma}\label{lemma-Heissubgroup}
A set $\Gamma=M\mZ^2\times m\mZ$, where $M$ is a matrix with integer entries, and $m>0$ is an integer, is a subgroup if and only if $m$ divides both entries of one of the rows of $M$.
\end{lemma}

Examples of group chains, where every subgroup is normal in $\cH$, and so by Theorem \ref{charact-groupchains} the action is homogeneous are given, for example, in \cite{LSU2014}. By Corollary \ref{cor-3-1} these actions have a trivial discriminant group.
 
We now give examples of the action of the Heisenberg group with non-trivial discriminant.

\begin{ex} \label{eq-wrHeisenberg}
{\rm 
Let $M_n=  \left( \begin{array}{cc} qp^n & pq^n \\ p^{n+1} & q^{n+1}\end{array}\right)$,  where $p, q \geq 2$ are distinct primes, and consider the action represented by a group chain $\G_0 = \cH$, $\{\G_n\}_{n \geq 1} = \{M_n \mZ^2 \times p\mZ\}_{n \geq 1}$.

\begin{lemma}
The discriminant group $\cD_x$ associated to the action corresponding to the chain $\{\G_n\}_{n \geq 0}$ is finite. The discriminant group $\cD_x$ is unstable.
\end{lemma}
 
 \proof We note that  for $x = (e\G_n)$ we have ${\ds K_x = \cap_n \G_n = \{0\} \times \{0\} \times p\mZ}$, which is a normal subgroup of $\cH$. Then by Proposition~\ref{prop-unstable} if $\cD_x$ is non-trivial, then it is unstable.
 
 By \cite{LSU2014}, a subgroup of $\cH$ which has the form $L\mZ^2 \times m\mZ$ is normal,  if $m$ divides every entry in $L$. So let $L_n=  \left( \begin{array}{cc} qp^n & p^2q^n \\ p^{n+1} & pq^{n+1}\end{array}\right)$, and consider a normal subgroup $L_n = B_n \mZ^2 \times p\mZ$, where $p,q$ are the same primes as in $M_n$. Then $L_n$ is a normal subgroup of $\G_n$ of index $p$. Since $p$ is a prime, $L_n$ is maximal in $\G_n$. Therefore, for every $n \geq 0$ the cardinality of $\G_n/L_n$ is $p$, and so the discriminant group $\cD_x$ has cardinality at most $p$ and is finite. Since the bonding maps $\G_n/L_n \to \G_{n-1}/L_{n-1}$ are clearly surjective, then $\cD_x$ is non-trivial.
 \endproof
 
 }
\end{ex}

 \begin{remark}
 {\rm  The   group chain in Example~\ref{eq-wrHeisenberg} is   weakly normal; see Example 5.10 in the  thesis by Dyer \cite{Dyer2015}.
 }
 \end{remark}

 \begin{ex}\label{main-6}
 {\rm Let $A_n =  \left( \begin{array}{cc} p^n & 0 \\ 0 & q^n\end{array}\right)$, where $p$ and $q$ are distinct primes, and consider the action represented by a group chain $G_0 =\cH $, $\{G_n\}_{n \geq 1} = \{A_n \mZ^2 \times p^n \mZ\}_{n \geq 1}$.

\begin{lemma}
The discriminant group $\cD_x$ associated to the action corresponding to the chain $\{G_n\}_{n \geq 0}$ is infinite. The discriminant group $\cD_x$ is unstable.\end{lemma}
 
 \proof  We note that for $x = (eG_n)$ the kernel $K _x= \bigcap_n G_n = \{e\}$ is trivial. Then by Proposition~\ref{prop-unstable} if $\cD_x$ is non-trivial, then it is unstable.

 By \cite{LSU2014}, a subgroup of $\cH$ which has the form $L\mZ^2 \times m\mZ$ is normal, if $m$ divides every entry in $L$. So let $E_n=  \left( \begin{array}{cc} p^n & 0 \\ 0 & p^n q^n\end{array}\right)$, and consider a subgroup 
  $$P_n = E_n \mZ^2 \times p^n\mZ = p^n \mZ \times p^nq^n \mZ \times p^n \mZ,$$ 
  where $p,q$ are the same primes as in $A_n$.  Then $P_n$ is a normal subgroup of $G_n$ of index $p^n$.  Since $P_n$ is contained in the maximal normal subgroup of $G_n$, and $p$ is a prime,  $P_n$ is maximal in $G_n$.  Therefore, for every $n \geq 0$ the cardinality of $G_n/P_n$ is $p^n$. Since the projections $G_{n+1}/P_{n+1}$ are clearly surjective,  then the discriminant group $\cD_x$ is infinite. 
 \endproof

 }
\end{ex}

\begin{remark}
 {\rm The   group chain in Example~\ref{main-6} is not weakly normal; see Example 5.14 in the  thesis by Dyer \cite{Dyer2015}. }
 \end{remark}

\subsection{Dihedral and the generalized dihedral groups}

We now give examples of actions of the dihedral group, which have infinite discriminant groups. 

\begin{ex}\label{example63}
{\rm
Let $\G = \mZ^2 = \{(a,b) \mid a,b \in \mZ\}$. Choose two distinct primes $p,q \in \mZ$, and define for $i \geq 1$
  \begin{align*}
  \G_i = \{(ap^{i}, bq^{i} ) \mid a,b, \in \mZ\} \ , \quad \textrm{and}\quad \G_i^T = \{(aq^{i}, bp^{i} ) \mid a,b, \in \mZ\} \ .
  \end{align*}
Let $H = \mZ_2 = \{1,t\}$, and define a homomorphism $\theta \colon   H \to {Aut}(\mZ^2)$ by
  $$\theta(t)(a,b) = (b,a).$$
Let $G = \G \rtimes H$, and $G_i = \G_i \times \{1\}$.  Given $g = ((a,b),t)$, we have $gG_ig^{-1} = \G_i^T \times \{1\}$,  and it follows that the core
  $$C_i = {\rm core}_G  \, G_i = G_i \cap \left(\G_i^T \times \{1\} \right) = \{(ap^iq^i,bp^iq^i,1) \mid a,b \in \mZ\}.$$
That is, $(aq^i,bp^i,1)$ and $(c q^i,dp^i,1)$ in $G_i$ are equal modulo $C_i$ if and only if $a = c \mod p^i$  and $b = d \mod q^i$. So $card(G_i/C_i) = p^iq^i$.

Now consider $\ds \cD_x = \lim_{\longleftarrow} \,\left\{\delta^i_{i-1} \colon   G_i/C_i \to G_{i-1}/C_{i-1} \right\}$. It is straightforward to that the bonding maps $\delta^i_{i-1} \colon  G_i/C_i \to G_{i-1}/C_{i-1}$ are surjective. Since $card(G_i/C_i)$ tends to infinity as $i \to \infty$, the discriminant group $\cD_x$ is infinite. Since $\bigcap_i \G_i$ is trivial, then the discriminant group $\cD_x$ is unstable.

Note that $\G \times \{1\}$ is an abelian subgroup of $G$, and $G_i \subset \G\times \{1\}$ for $i \geq 1$. Then $G_i$ is normal in $G_1$ for $i \geq 1$, and action is weakly normal.
}
\end{ex}

\begin{ex}\label{example63-1}
{\rm
Example \ref{example63} can be generalized to produce a family of examples as follows. Let $\G =\mZ^n$, and let $H = A_n$, the permutation group on $n$ elements. Set $G = \G \rtimes H$. Now choose $n$ distinct primes $p_1,\ldots,p_n$, and define $\G_i = \{(a_1p_1^i, a_2 p_2^i,\ldots,a_np_n^i)\}$ and $G = \G \times \{1\}$. A computation, similar to the one in Example \ref{example63} shows that the equicontinuous system given by the action of $G$ on the inverse limit $G_{\infty}$, associated to the group chain $\{G_i\}_{i \geq 0}$, is infinite. Since $\bigcap_i \G_i$ is trivial, then the discriminant group $\cD_x$ is unstable. Since $\G \times \{1\}$ is an abelian subgroup of $G$, and $G_i \subset \G \times \{1\}$ for $i \geq 1$, the action is weakly normal.

}
\end{ex}

 \vfill
 \eject
 
\appendix

\section{Group chains and inverse limit representations}\label{appendix}

In this section, given a minimal equicontinuous action $(X,G,\Phi)$, we construct the  associated group chain $\{G_i\}_{i \geq 0}$ using the methods of \cite{ClarkHurder2013}.  The group chain $\{G_i\}_{i \geq 0}$ so constructed is dependent on the choice of a collection of clopen partitions of $X$ and on the choice of a point in $X$. 

Let $d$ be a metric on $X$. Recall that the action $(X,G,\Phi)$ is equicontinuous if for any $\e>0$ there exists $\delta >0$ such that for any $g \in G$ and any $x,y \in X$ such that $d(x,y) < \delta$ we have $d(g \cdot x,g \cdot y)) < \e$. The constant $\delta$ is called an \emph{equicontinuity constant} for the $\e>0$. Recall that $(X,G,\Phi)$ is minimal if the orbit of any $x \in X$ under the action of $G$ is dense in $X$.

\begin{prop}\label{prop-AFpres-1}
Let $(X,G,\Phi)$ be a minimal equicontinuous action, and let $\whx \in X$ be a point. Then there exists an infinite sequence of closed and open sets $X = V_1^{(0)} \supset V_1^{(1)} \supset V_1^{(2)}  \supset \cdots$ with $ \ds \{\whx\} = \bigcap_{i \geq 1} V_1^{(i)} $, which have the following properties.
\begin{enumerate}
\item For $i \geq 1$ the collection $\cP_i =\{g \cdot V_1^{(i)} \mid g \in G\}$ is a finite partition of $X$ into clopen sets. \item We have ${ diam}(V_k^{(i)} ) \leq 2^{-i}$ for all $1 \leq k \leq \kappa_i$.
\item The collection of elements which fix $V_1^{(i)}$, that is,
 $$G_i = \{g \in G \mid g \cdot V_1^{(i)} = V_1^{(i)}\},$$
 is a subgroup of finite index in $G$. More precisely, $|G:G_i| = { card}(\cP_i)$.
\item There is a homeomorphism $\ds \phi \colon  X \to   G_{\infty} = \lim_{\longleftarrow} \,\left\{G/G_{i+1} \to G/G_i \right\}$ equivariant with respect to the action of $G$ on $X$ and $G_{\infty}$, which maps $\whx$ onto $(eG_i) \in G_{\infty}$.
\end{enumerate}
\end{prop}

\proof Let $0 < \e_1 < 1/2$, and $0 < \e_{i+1} <  \e_{i}/2$ for $i \geq 1$, be a sequence of real numbers. Recall that for any two sets $U,V \subset X$ the distance between $V$ and $U$ is given by
  \begin{align*} { dist}(U,V) = \inf\{d(u,v) ~|~ u \in U, v \in V~\}.\end{align*}

Choose a finite clopen partition of $X$, $\cW^{(1)}=\{W_1^{(1)},...,W_n^{(1)} \}$ such that ${ diam}(W_i^{(1)})<\e_1$. Since the sets in $\cW^{(1)}$ are disjoint, and $\cW^{(1)}$ is finite, there exists a number $\tilde{\e}_1>0$ such that
  \begin{align} \label{vetilde}\tilde{\e}_1 < \min \{ \e_1, \{ { dist} (W_i^{(1)}, W_j^{(1)}) \mid i \neq j , 1 \leq i,j \leq n\} \}.\end{align}
Let $\delta(\tilde{\e}_1)$ be an equicontinuity constant for $\tilde{\e}_1$. 
For each $x \in X$, we define a coding function $C_x^{(1)} \colon  G \to \{1,...,n \}$ as follows. 
 \begin{align} \label{coding-function} C_x^{(1)} (g) = i ~ \textrm{ if  and only if }g \cdot x \in W_i^{(1)}.\end{align}
That is, the coding function of $x$ tells us which set of the partition $\cW^{(1)}$ contains $g \cdot x$.

For each $x \in X$, consider a set of points in $X$ which have the same coding functions as $x$, that is,
 \begin{align}\label{eq-Vx}V_x=\{ y \in X \mid  C_y^{(1)}=C_x^{(1)} \}.\end{align}
If $x \in W_i^{(1)}$, then $C_x^{(1)}(e) = i$, where $e \in G$ is the identity element. Then $V_x \subseteq W_i^{(1)}$, and so the collection $\{V_x\}_{x \in X}$ is a refinement of $\cW^1$. Then by \eqref{vetilde} we have $diam(V_x) < \e_1$.

\begin{lemma}\label{clopensets}
The collection $\cP_1$ is a finite clopen partition of $X$. 
\end{lemma}

\proof First we show that for any $x \in X$, the set $V_x$ is open, that is, for any $x'$ with $d(x,x')<\delta(\tilde{\e}_1)$, we have $$C_x^{(1)}(g)=C_{x'}^{(1)}(g) \textrm{ for all } g \in G.$$ 
Indeed, suppose $g \in G$, and let
 $d(x,x')<\delta(\tilde{\e}_1)$. Since the action is equicontinuous, $d(g \cdot x, g \cdot x')<\tilde{\e}_1$. Then if $g \cdot x \in W^{(1)}_i$, we must have $g \cdot x' \in W_i^{(1)}$ as well, since by \eqref{vetilde} $\tilde{\e}_1$ is less than the distance between any two sets in $\cW^{(1)}$. Thus  $C_x^{(1)}(g)=C_{x'}^{(1)}(g)=i$. 

We now show that the sets $\{V_x\}$ are disjoint, that is, if $V_x \cap V_y \neq \emptyset$, then $V_x=V_y$. Let $z \in V_x \cap V_y$. Since $z \in V_x$, then for any $g \in G$ we have
   $$C_x^{(1)}(g)=C_z^{(1)}(g).$$
Since $z \in V_y$, then for any $g \in G$ we have
  $$C_y^{(1)}(g)=C_z^{(1)}(g).$$
Therefore, $C_x^{(1)}(g)=C_y^{(1)}(g)$, which implies that $V_x = V_y$. Since every $V_x$ is an open set, the collection $\cP_1$ is an open cover of $X$ by disjoint sets. Since $X$ is compact, this cover is finite. Then the complement of each $V_x$ is a finite union of open sets, and so must be open. Therefore, $V_x$ is also closed.
\endproof

Let $\kappa_1$ be the cardinality of $\cP_1$, and let $\cP_1 = \{V_1^{(1)},\ldots,V_{\kappa_1}^{(1)}\}$ with $\whx \in V_1^{(1)}$. 

\begin{lemma} \label{permutation lemma} The action of $G$ permutes the sets of $\cP_1$, that is, for each $V_i^{(1)} \in \cP_1$ and each $g \in G$, $g \cdot V_i^{(1)}=V_j^{(1)}$ for some $1\leq j \leq \kappa_1$. 
\end{lemma}

\proof
Suppose $x,y \in V_i^{(1)}$. By the definition of $V_i^{(1)}$, $g \cdot x, g \cdot y $ are in the same $W_k$ for all $g \in G$. 
Let $g \cdot x\in V_j^{(1)} \subset W_k$ and $g \cdot y \in V^{(1)}_m \subset W_k$ for some $g \in G$, and suppose $j \ne m$. Then for some $g' \in G$
   $$C_{g \cdot x}(g') = \ell, \textrm{ and } C_{g \cdot y}(g') = \ell',$$ 
 where $\ell \ne \ell'$. Then $g'g \cdot x \in W_\ell$ and $g' g \cdot  y \in W_{\ell'}$, where $W_\ell$ and $W_{\ell'}$ are distinct sets in $\cW^{(1)}$. This is not possible, since $x,y \in V_i^{(1)}$, the same set of the partition $\cP_1$. It follows that $j = m$.
 \endproof

Let $G_1 = \{ g \in G \mid g \cdot V_1^{(1)} = V_1^{(1)}\}$ be the isotropy subgroup of the action at $V_1^{(1)}$, that is, $G_1$ is the set of elements which fix $V_1^{(1)}$. Lemma \ref{lemma-cosets} shows that $G_1$ has finite index in $G$.

\begin{lemma}\label{lemma-cosets}
There is a natural bijective map $\tilde{\phi}_1 \colon  \cP_1 \to G/G_1$, and so
  $$|G:G_1| = { card}(\cP_1).$$ 
\end{lemma}

\proof Define $\tilde{\phi}_1 \colon \cP_1 \to G/G_1$ by $\tilde{\phi}_1(V_k^{(1)}) = gG_1$ if and only if $g \cdot V_1^{(1)} = V^{(1)}_k$. To see that this map is well-defined, suppose  $h \in gG_1$. Then $g^{-1}h \in G_1$, and so $g^{-1}h \cdot V_1^{(1)} = V_1^{(1)}$. Then $h \cdot V_1^{(1)} =g(g^{-1}h) \cdot V_1^{(1)}= g \cdot V_1^{(1)}$. On the other hand, let $g \cdot V_1^{(1)} = h \cdot V_1^{(1)}$ for $g,h \in G$. Then $g^{-1}h \in G_1$, and so $g \in hG_i$.
\endproof

We now perform the inductive step of the construction. Suppose we have a group chain 
  $G=G_0 \supset G_1 \supset \cdots \supset G_i$ and a sequence of refining partitions $\cP_1, \cP_2, \ldots, \cP_i$ satisfying the conditions of Proposition \ref{prop-AFpres-1}, that is, ${ diam}(V_k^{(i)}) < \e_i$ for each $V_k^{(i)} \in \cP_i$, the group $G_i$ is the isotropy group at $V_1^{(i)} \in \cP_i$, $|G:G_i| = { card} ( \cP_i)$, and the sets in $\cP_i$ are permuted by the action of $G$. 

Let $\cW^{(i+1)}$ be a partition of $X$ refining $\cP_i$, such that ${ diam}(W_k) < \e_{i+1}$ for all $W_k \in \cW^{(i+1)}$. For each $x \in X$, define the coding functions $C_x^{(i+1)} \colon  G_{i} \to \cW^{(i+1)}$ by
\begin{align*} C_x^{(i+1)}(g) = k \textrm{ if and only if }g \cdot x \in W_k,\end{align*}
and let $V_x = \{ y \in X \mid C^{(i+1)}_x (g)= C^{(i+1)}_y(g) \textrm{ for every }g \in G\}$. 
By an argument, similar to the one in Lemmas \ref{clopensets}, \ref{permutation lemma} and \ref{lemma-cosets}, one shows that $\cP_{i+1}= \{V_x\}_{x \in X}= \{V_1^{(i+1)}, \ldots, V_{\kappa_i}^{(i+1)}\}$ is a clopen partition of $X$ with $\whx \in V_1^{(i+1)}$ and the following properties. 
\begin{enumerate}
\item ${ diam}(V_k^{(i+1)}) < \e_{i+1}$ for every $V_k^{(i+1)} \in \cP_{i+1}$,
\item The partition $\cP_{i+1}$ refines $\cW^{(i+1)}$ and so refines $\cP_i$,
\item The action of $G$ permutes the sets of $\cP_{i+1}$, 
\item $G_{i+1} = \{g \in G_i \mid g \cdot V_1^{(i+1)} = V_1^{(i+1)}\}$ is a subgroup of $G_i$.
\item There is a bijective map $\tilde{\phi}_{i+1} \colon  \cP_i \to G/G_i$, and so $|G: G_{i+1}| = { card} (\cP_{i+1})$ and $G_{i+1}$ has finite index in $G$.
\end{enumerate}

Now suppose we have a collection $\{\cP_i\}_{i \geq 0}$ of finite clopen partitions of $X$ as above, and 
let $\theta^i_{i-1}\colon G/G_i \to G/G_{i-1}$ be the inclusion of cosets. Then the inverse limit 
  $$G_{\infty} =\lim_{\longleftarrow} \, \left\{ \theta^i_{i-1} \colon  G/G_i \to G/G_{i-1} \right\} \subset \prod_i G/G_i$$ 
  is a Cantor set in the relative topology from the product topology on $\prod_i G/G_i$. The left action of $G$ on $G_{\infty}$ is defined by
  \begin{align*} G \times G_{\infty} \to G_{\infty} : (h, (g_iG_i)) \mapsto (hg_iG_i). \end{align*}
For every $\cP_i$, let $\jmath_i (x) = V_k^{(i)}$ if and only if $x \in V_k^{(i)}$. Then there are $G$-equivariant maps
  \begin{align} \label{eq-finiterep}
   \phi_i = \tilde{\phi}_i \circ \jmath_i  \colon   X \to G/G_i,
   \end{align}
where $\tilde{\phi}_i$ is defined in Lemma \ref{lemma-cosets}. It is readily verified that the mappings $\phi_i$ are compatible with the bonding maps $\theta^i_{i-1}$, and so there is a continuous map
  \begin{align}\label{map-phi}
  \phi \colon  X \to G_{\infty} = \lim_{\longleftarrow} \,\left\{\theta_i \colon  G/G_i \to G/G_{i-1} \right\}.
  \end{align}
Since every coset space $G/G_i$ is finite, and $\phi_i$ are surjective maps, by \cite[Corollary 1.1.6]{RZ2000} the map $\phi$ is surjective. 

Let $x,y \in X$, and choose $i \geq 1$ such that $\e_i < d(x,y)$. Since sets in $\cP_i$ have diameter less than $\e_i$, $x$ and $y$ are in different sets of $\cP_i$, and so $\phi_i(x) \ne \phi_i(y)$, and $\phi(x) \ne \phi(y)$. Thus $\phi$ is injective, and in fact a homeomorphism.
Since the action of $G$ permutes sets in $\cP_i$, the mapping $\phi \colon  X \to G_{\infty}$ is equivariant with respect to the action of $G$.
\endproof

Given an equicontinuous minimal group action $(X,G,\Phi)$ on a Cantor set $X$, Proposition \ref{prop-AFpres-1} associates to it a nested group chain $\{G_i\}_{i \geq 0}$. 
This group chain depends on the choice of the collection of refining clopen partitions $\{\cP_i\}_{i \geq 0}$, and on the choice of a point $\whx \in X$.


\end{document}